\documentclass{article}
\usepackage{float}
\usepackage{graphicx}
\usepackage{enumitem}
\usepackage{amsmath,amssymb,amsthm,bm,mathtools}
\usepackage{hyperref}
\usepackage{cite}
\graphicspath{{figures/}}
\hypersetup{
    bookmarks=true,         
    unicode=false,          
    pdftoolbar=true,        
    pdfmenubar=true,        
    pdffitwindow=false,     
    pdfstartview={FitH},    
    pdfnewwindow=true,      
    colorlinks=true,        
    linkcolor=blue,         
    citecolor=blue,         
    filecolor=magenta,      
    urlcolor=blue           
}
\newcommand{\dsp}{\displaystyle}
\newcommand{\Zbb}{\mathbb{Z}}
\newcommand{\Qbb}{\mathbb{Q}}
\newcommand{\Rbb}{\mathbb{R}}
\newcommand{\Cbb}{\mathbb{C}}
\newcommand{\ii}{\mathrm{i}}
\newcommand{\mcz}{\mathcal{Z}}
\newcommand{\mcm}{\mathcal{M}}

\DeclareMathOperator{\arcsinh}{arcsinh}

\DeclareMathOperator{\Tr}{Tr}
\newtheorem{theorem}{Theorem}[section]
\newtheorem{lemma}[theorem]{Lemma}

\newtheorem{remark}[theorem]{Remark}
\begin{document}
\begin{center}
{\Large \bf Physical configurations of a cell doublet with line tension, a theoretical study.}\\
\bigskip
\bigskip
\textit{Fabrice Delbary}\\
\bigskip
\bigskip
Collège de France, 11 Place Marcelin-Berthelot, 75231 Paris Cedex 05, France
\end{center}
\bigskip
\bigskip
\bigskip
\bigskip
\section{Introduction}
Although being a very simple configuration, a cell doublet can appear to be useful to perform simple preliminary tests of numerical algorithms, both in the framework of forward and inverse problems. For a cell doublet with surface tensions only, the link between the geometry and the surface tensions is well-known \cite{barrat,chaikin_lubensky,gennes_brochard-wyart_quere,israelachvili,safran} and the double bubble problem for immiscible fluids has been mathematically proven, first for certain volumes and nearly-unit surface tensions \cite{hutchings_morgan_ritore_ros} and later in the general case \cite{lawlor}. In the first part of the manuscript, in a simpler framework where it is \textit{a priori} assumed that the interfaces of the cell doublet are spherical caps, we derive the equations giving the geometrical parameters for known surface tensions and volumes. Considering the proof in \cite{lawlor}, we obviously recover the fact that the configuration is unique and the Young-Laplace and Young-Dupré laws. The solution of the geometrical parameters can then be written in a semi-explicit way. The main purpose here is to introduce the method we use in order to write the configuration of a cell doublet with line tension in a semi-explicit way as well. In the second part, we show how to compute the geometrical parameters of a cell doublet with surface tensions when the pressures instead of the volumes are known. This case does not present any difficulty since the Young-Laplace and Young-Dupré laws immediately give simple relations between the geometrical parameters and the tensions and pressures. In the third part, we consider a cell doublet with surface tensions and a line tension. Using the same method as for the cell doublet without line tension, it seems that the configuration is unique. Moreover, the geometrical parameters can be written in a semi-explicit way as well.
\section{Description}
We consider two cells $1$ and $2$ joined by an interface $3$ and want to study the configuration of the system at minimum energy. The cell surfaces are modeled as three spherical caps $\mathcal{S}_1$, $\mathcal{S}_2$, $\mathcal{S}_3$ with respective centers $C_1$, $C_2$, $C_3$, where $\mathcal{S}_3$ defines the interface $3$ and where the surface of cell $k=1,2$ is given by $\mathcal{S}_k\cup\mathcal{S}_3$. We consider the oriented $(Ox)$ axis whose origin $O$ is the center of the base circle of $S_3$ and whose direction is given by the vector ${\bf e_x}=(C_2-C_1)/|C_2-C_1|$. We denote by $\mathcal{C}_0$ (respectively $\mathcal{D}_0$), the unit circle (respectively unit disk) of center $O$ orthogonal to ${\bf e_x}$.\\
\\
For $x\in\Rbb^*$ and $h\in\Rbb$, we denote by $S_{x,h}$ the spherical cap with base circle $\dsp h\mathcal{C}_0$ and apex $x{\bf e_x}$. $S_{0,h}$ is defined as the disk $\dsp h\mathcal{D}_0$. The surfaces $1,2,3$ are then given by $S_{x_1,h}$, $S_{x_2,h}$, $S_{x_3,h}$ for $x_1,x_2,x_3,h\in\Rbb$ with $x_1<x_3<x_2$.
\begin{remark}
For $h\in\Rbb$, $\pm h$ define the same spherical caps, hence the same configuration of the cell doublet, but it is more convenient to keep it this way than restraining $h\in\Rbb^+$.
\end{remark}
One might want to introduce the radii of the spherical caps as well as their aperture in order to express in a more common way the force balance. Hence, we define for $k\in\{1,2,3\}$
\begin{equation}\label{zar}
z_k=\frac{x_k}{|h|}\quad,\quad \alpha_k=2\arctan z_k\quad,\quad r_k=\frac{|h|}{\sin\alpha_k}\quad\mbox{(signed radius)}.
\end{equation}
For sake of simplicity, the cosines and sines are abbreviated for $k\in\{1,2,3\}$ as
\begin{equation}
c_k=\cos\alpha_k\qquad,\qquad s_k=\sin\alpha_k.
\end{equation}
We then have for $k\in\{1,2,3\}$
\begin{equation}
c_k=\frac{1-z_k^2}{1+z_k^2}\qquad,\qquad s_k=\frac{2z_k}{1+z_k^2}.
\end{equation}
The tension of each surface $1,2,3$ is respectively denoted by $t_1,t_2,t_3\in\Rbb_+^*$. Note that the angles of tension forces with respect to the $(Ox)$ axis are respectively given by $\alpha_k-\pi/2$, $k\in\{1,2,3\}$. The volumes of cells $1$ and $2$ are respectively given by
\begin{align}
V_1&=\frac{\pi}{6}\left(x_3(3h^2+x_3^2)-x_1(3h^2+x_1^2)\right),\\
V_2&=\frac{\pi}{6}\left(x_2(3h^2+x_2^2)-x_3(3h^2+x_3^2)\right),
\end{align}
Which can also be written with the variables $z$ as
\begin{align}
V_1&=\frac{\pi}{6}h^3\left(z_3(z_3^2+3)-z_1(z_1^2+3))\right),\\
V_2&=\frac{\pi}{6}h^3\left(z_2(z_2^2+3)-z_3(z_3^2+3)\right),
\end{align}
Or with the variables $\alpha$ as
\begin{align}
V_1&=\frac{\pi}{3}h^3\left(\frac{s_3(2+c_3)}{(1+c_3)^2}-\frac{s_1(2+c_1)}{(1+c_1)^2}\right),\\
V_2&=\frac{\pi}{3}h^3\left(\frac{s_2(2+c_2)}{(1+c_2)^2}-\frac{s_3(2+c_3)}{(1+c_3)^2}\right).
\end{align}
We define the difference of angles by
\begin{equation}
\varphi_k=\alpha_{k+1}-\alpha_{k-1},\quad\mbox{for }k\in\{1,2,3\},
\end{equation}
Where indices are considered modulo $3$ (in $\Zbb/3\Zbb$), that is $\alpha_0=\alpha_3,\alpha_4=\alpha_1$, etc.. These angles satisfy
\begin{equation}
\varphi_1+\varphi_2+\varphi_3=0.
\end{equation}
\begin{figure}[H]
\begin{center}
\includegraphics[scale=0.75]{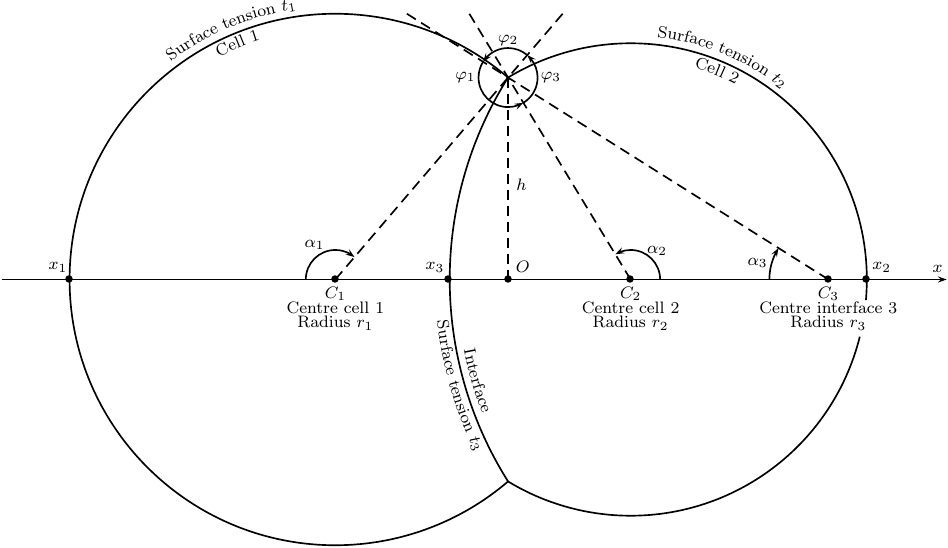}
\caption{Cell doublet configuration. Angles are taken in $(-\pi,\pi)$.}\label{cellconf}
\end{center}
\end{figure}
\section{Configuration at minimum energy under volume constraints}\label{cellconfminen}
For constrained volumes, the pressure forces exert no work under a change of configuration and the energy $E$ of the system is simply given by the sum of surface tensions times the spherical caps areas, that is
\begin{equation}
E=\pi(t_1(x_1^2+h^2)+t_2(x_2^2+h^2)+t_3(x_3^2+h^2))
\end{equation}
That is
\begin{equation}
E=\pi(t_1x_1^2+t_2x_2^2+t_3x_3^2+t_{\mathrm{s}} h^2),
\end{equation}
Where
\begin{equation}
t_{\mathrm{s}}=t_1+t_2+t_3.
\end{equation}
Note that the energy can also be written using the $z$ or $\alpha$ variables and we have
\begin{equation}
E=\pi h^2(t_1z_1^2+t_2z_2^2+t_3z_3^2+t_{\mathrm{s}}),
\end{equation}
\begin{equation}
E=2\pi h^2\left(\frac{t_1}{1+c_1}+\frac{t_2}{1+c_2}+\frac{t_3}{1+c_3}\right)
\end{equation}
\subsection{Equations for the critical points of the constrained energy}\label{sysexind}
Since the problem of finding the critical points of the energy is purely algebraic (all the equations are given by polynomials), although not mandatory, it is interesting to describe it in an algebraic way. You can refer to \cite{cox_sturmfels_manocha_sederberg_kramer_laubenbaches_thomas,michalek_sturmfels,sturmfels_1995,sturmfels_2002}, but as a quick reminder, and roughly speaking, an affine variety is an implicit manifold whose defining equations are polynomials. It can happen that the affine variety has singularities (it is not a manifold at these points), the dimension of an affine variety is however defined as its maximal dimension at regular points when viewed as a manifold. To end this very short, incomplete and non exhaustive presentation, an ideal is a set of polynomials which is a group for the sum and such that the multiplication of a polynomial in the ideal by any other polynomial (not necessary in the ideal) is sill in the ideal. Hence, zeros of a polynomial system can also be characterized as the zeros of the ideal generated by the polynomials of the system (smallest ideal containing the polynomials).\\
\\
We consider the problem of finding the minimizer of the energy $E$ (or equivalently $E/\pi$) under volume constraints $V_1(x_1,x_2,x_3,h)=v_1$ and $V_2(x_1,x_2,x_3,h)=v_2$, where $v_1,v_2\in\Rbb_+^*$, that is, we want to find the minimizer of the polynomial $E(X_1,X_2,X_3,Y)$ on the affine variety $\mcm$ of dimension $2$ of $\Rbb^4$ defined by the ideal $I_{\mcm}=\langle Q_1,Q_2\rangle$ where $Q_1,Q_2\in\Rbb[X_1,X_2,X_3,Y]$ are given by
\begin{align}
Q_1&=-X_1(3Y^2+X_1^2)+X_3(3Y^2+X_3^2)-w_1,\\
Q_2&=X_2(3Y^2+X_2^2)-X_3(3Y^2+X_3^2)-w_2,
\end{align}
Where for $k\in \{1,2\}$, $w_k=6v_k/\pi$. We also denote $v_3=v_1+v_2$ and $w_3=6v_3/\pi$. Note that $x_1<x_2<x_3$ is implicitly imposed by the volume constraints.
\begin{remark}
An homothety of the surface tensions does not change the configuration of the system at minimum energy. Similarly, an homothety of factor $\lambda^3$ of the volumes $v_1,v_2$ simply scales the spatial coordinates of the system at minimum energy by a factor $\lambda$, thus, the angles $\alpha_1,\alpha_2,\alpha_3$ remain unchanged.
\end{remark}
\begin{remark}
It is an easy task to solve numerically the problem and find the configuration of the double cell. However, we would like to have explicit formulae that could be used directly. Since the problem has been written in terms of a minimization of a polynomial on an affine variety, the solution will hence be given as the solution to polynomial system. Using a Gröbner basis \cite{sturmfels_1995} (computed for instance with a Computer Algebra System like Singular \cite{singular}), the solution is then simply given by taking the successive roots of polynomials whose coefficients depend on the previous parameters given by the previous polynomials. However, note that although semi-explicit, this way of solving the system is not the most appropriate for numerical applications since Gröbner basis often give large polynomial coefficients which can dramatically destroy the precision. Nonetheless, in our simple case of a double cell, we can easily find an univariate polynomial whose only real root is one of the parameter $z_1,z_2,z_3$. The other $z$ parameters are then given by Möbius transformations depending on the surface tensions $t_1,t_2,t_3$ and $h$ can be computed using the volumes and the $z$ parameters.
\end{remark}
\begin{remark}
For $q\in\Rbb$, the polynomial $X(X^2+3)-2q$ has a unique real root $z$ given by
\begin{equation}
z=\sqrt[3]{q+\sqrt{1+q^2}}-\frac{1}{\sqrt[3]{q+\sqrt{1+q^2}}}=\sqrt[3]{q+\sqrt{1+q^2}}+\sqrt[3]{q-\sqrt{1+q^2}},
\end{equation}
Which can also be written
\begin{equation}
z=\mcz(q),
\end{equation}
Where
\begin{equation}
\mcz(q)=2\sinh\left(\frac{1}{3}\arcsinh(q)\right).
\end{equation}
Hence, $\mcm$ can be easily parameterized using $h$ and one of the other three variables $x_1,x_2,x_3$, indeed, $\mcm$ is any of the set of points $x$ respectively defined by
\begin{subequations}
\begin{align}\label{paramm}
x=\left(x_1,h\mcz\left(\frac{x_1(3h^2+x_1^2)+w_3}{2h^3}\right),h\mcz\left(\frac{x_1(3h^2+x_1^2)+w_1}{2h^3}\right),h\right),\\
x=\left(h\mcz\left(\frac{x_2(3h^2+x_2^2)-w_3}{2h^3}\right),x_2,h\mcz\left(\frac{x_2(3h^2+x_2^2)-w_2}{2h^3}\right),h\right),\\
x=\left(h\mcz\left(\frac{x_3(3h^2+x_3^2)-w_1}{2h^3}\right),h\mcz\left(\frac{x_3(3h^2+x_3^2)+w_2}{2h^3}\right),x_3,h\right).
\end{align}
\end{subequations}
The limit is well-defined as $h\to0$, giving in fact parameterizations of the closed subvariety $\mcm_0$ of $\mcm$ of points with $h=0$ (algebraic curve of the hyperplane $\Rbb^3\times\{0\}$), that is, the points $(x_1,x_2,x_3,0)\in\Rbb^4$ such that
\begin{align*}
-x_1^3+x_3^3&=w_1,\\
x_2^3-x_3^3&=w_2.
\end{align*}
The respective parameterizations are given by
\begin{align}\label{paramm0}
\mcm_0&=\left\{\left(x_1,\sqrt[3]{x_1^3+w_3},\sqrt[3]{x_1^3+w_1},0\right)\,:\,x_1\in\Rbb\right\},\\
\mcm_0&=\left\{\left(\sqrt[3]{x_2^3-w_3},x_2,\sqrt[3]{x_2^3-w_2},0\right)\,:\,x_2\in\Rbb\right\},\\
\mcm_0&=\left\{\left(\sqrt[3]{x_3^3-w_1},\sqrt[3]{x_3^3+w_2},x_3,0\right)\,:\,x_3\in\Rbb\right\},
\end{align}
\end{remark}
Remark that $E$ is a positive definite quadratic form of $\Rbb^4$ and $\mcm$ is closed, hence the minimum $E$ on $\mcm$ is reached at some point. Moreover, it can be easily checked that the (real) singular locus $\mcm_{\mathrm{sing}}$ of $\mcm$ is empty, thus the minimum of $E$ on $\mcm$ is a critical point. The singular locus of $\mcm$ is the set points of $\mcm$ which are singular, that is points such that $\mcm$ is not a differentiable manifold at these locations, the tangent space has not dimension $2$ anymore, that is, if we denote by $\dsp Q=\begin{pmatrix}Q_1\\Q_2\end{pmatrix}$ and by $J_Q$ the Jacobian of $Q$ (basis of the orthogonal to the tangent plane), these points vanish the $2\times 2$ minors of $J_Q$ and of course the polynomials of $I_{\mcm}$ since the points have to be on the variety. The ideal $I_{\mcm_{\mathrm{sing}}}$ generated by these polynomials can then be written
\begin{equation}
I_{\mcm_{\mathrm{sing}}}=I_{\mcm}+\langle 2\times2\mbox{ minors of }J_Q\rangle,
\end{equation}
That is
\begin{equation}
\begin{aligned}
  I_{\mcm_{\mathrm{sing}}}&=I_{\mcm}+\left\langle Y(X_2-X_3)(Y^2+X_1^2),Y(X_3-X_1)(Y^2+X_2^2),\right.\\
  &\left.Y(X_1-X_2)(Y^2+X_3^2),(Y^2+X_2^2)(Y^2+X_3^2),\right.\\
  &\left.(Y^2+X_3^2)(Y^2+X_1^2),(Y^2+X_1^2)(Y^2+X_2^2)\right\rangle.
\end{aligned}
\end{equation}
We can see that the zeros of $\langle 2\times2\mbox{ minors of }J_Q\rangle$ are quadruples $(x_1,x_2,x_3,h)\in\Cbb^4$ such that $x_k=x_\ell=\pm\ii h$ for a couple $(k,\ell)\in\{1,2,3\}^2,k\neq\ell$, hence $\mcm_{\mathrm{sing}}=\emptyset$.\\
\\
We now want to find the critical points of the constrained energy. These points are characterized as those such that the gradient of the energy is orthogonal to tangent plane to $\mcm$, that is, the gradient of the energy lies in the space generated by $\nabla Q_1,\nabla Q_2$. This is exactly what is done when working with the Lagrangian $\mathcal{L}=E-P_1Q_1-P_2Q_2$. Deriving it, we search the points $(x_1,x_2,x_3,h)$ and coefficients $P_1,P_2$ (pressures) such that $\nabla E$ can be written as a linear combination of basis vectors in the orthogonal of the tangent plane, that is $\nabla E=P_1\nabla Q_1+P_2\nabla Q_2$.\\
\\
We can also proceed by searching for the points of $\mcm$ vanishing the $3\times 3$ minors of the gradient $J_{E_Q}$ of $\dsp E_Q=\begin{pmatrix}E\\Q\end{pmatrix}$, the critical ideal $\dsp I_{E,Q}^{\mathrm{crit}}$ of the volume constrained energy (polynomials whose zeros are the critical points of the constrained energy) is then given by
\begin{equation}
I_{E,Q}^{\mathrm{crit}}=I_{\mcm}+\langle 3\times3\mbox{ minors of }J_{E_Q}\rangle
\end{equation}
That is
\begin{equation}
\begin{aligned}
I_{E,Q}^{\mathrm{crit}}=&I_{\mcm}+\left\langle t_1X_1(Y^2+X_2^2)(Y^2+X_3^2)+t_2X_2(Y^2+X_3^2)(Y^2+X_1^2)+\right.\\
&\left.t_3X_3(Y^2+X_1^2)(Y^2+X_2^2),2t_2X_2Y(X_1-X_2)(Y^2+X_3^2)+\right.\\
&\left.2t_3X_3Y(X_1-X_3)(Y^2+X_2^2)+t_{\mathrm{s}}Y(Y^2+X_2^2)(Y^2+X_3^2),\right.\\
&\left.2t_3X_3Y(X_2-X_3)(Y^2+X_1^2)+2t_1X_1Y(X_2-X_1)(Y^2+X_3^2)+\right.\\
&\left.t_{\mathrm{s}}Y(Y^2+X_3^2)(Y^2+X_1^2),2t_1X_1Y(X_3-X_1)(Y^2+X_2^2)+\right.\\
&\left.2t_2X_2Y(X_3-X_2)(Y^2+X_1^2)+t_{\mathrm{s}}Y(Y^2+X_1^2)(Y^2+X_2^2)\right\rangle.
\end{aligned}
\end{equation}
Hence, the algebraic degree of the problem is $75$ (number of zeros of the critical ideal in $\Cbb^4$), however, most of critical points are without interest for our considerations since they are either complex or do not correspond to a minimum of the constrained energy but to a maximum or a saddle point. But we can already see that for $h\neq0$, any critical point of the constrained energy has to verify the volume constraints $Q_1,Q_2$ and
\begin{align}\label{young_dupre}
t_1\frac{h^2-x_1^2}{h^2+x_1^2}+t_2\frac{h^2-x_2^2}{h^2+x_2^2}+t_3\frac{h^2-x_3^2}{h^2+x_3^2}&=0,\\
t_1\frac{x_1}{h^2+x_1^2}+t_2\frac{x_2}{h^2+x_2^2}+t_3\frac{x_3}{h^2+x_3^2}&=0,
\end{align}
That is
\begin{align}\label{young_dupre_angles}
t_1c_1+t_2c_2+t_3c_3&=0,\\
t_1s_1+t_2s_2+t_3s_3&=0,
\end{align}
Which is the force balance on the junction circle (Young-Dupré law) \cite{barrat,chaikin_lubensky,gennes_brochard-wyart_quere,israelachvili,safran}. It shows in particular that $x_1$, $x_2$, $x_3$ can not have the same sign (no non-physical moon-shape configuration for the system). Geometrically speaking, since the sum of tension forces is $0$, they form a triangle whose lengths are $(t_1,t_2,t_3)$. This implies that $(\varphi_1,\varphi_2,\varphi_3)$ is uniquely determined up to a minus sign. In addition, since we have $x_1<x_3<x_2$ because of volume constraints (cell $1$ on the left and cell $2$ on the right), we then have $x_1<0$ and $x_2>0$ which determines the sign of $(\varphi_1,\varphi_2,\varphi_3)$. Moreover, as we will detail later, the tensions have to fulfill the triangle inequalities.
\begin{figure}[H]
\centering
\includegraphics[width=5cm]{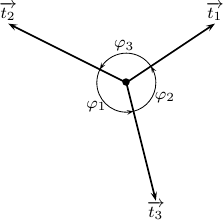}\hspace{2cm}
\includegraphics[width=5cm]{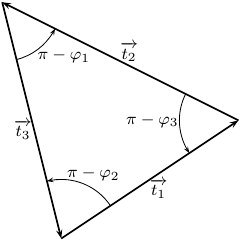}
\caption{Tension forces and angles}\label{forces_triangle}
\end{figure}
Still considering a possible critical point with $h\neq0$, the triangle formed by the tension forces has angles $\pi-\varphi_k$, $k\in\{1,2,3\}$ and the law of cosines (dot product of each couple of forces) gives for $k\in\{1,2,3\}$
\begin{align}
\cos\varphi_k&=-\frac{t_{k+1}^2+t_{k-1}^2-t_k^2}{2t_{k+1}t_{k-1}},\label{coslaw}\\
&=1-\frac{(t_{k+1}+t_{k-1}+t_k)(t_{k+1}+t_{k-1}-t_k)}{2t_{k+1}t_{k-1}},\nonumber\\
&=-1+\frac{(t_k+t_{k+1}-t_{k-1})(t_k+t_{k-1}-t_{k+1})}{2t_{k+1}t_{k-1}},\nonumber\\
\cos\left(\frac{\varphi_k}{2}\right)&=\sqrt{\frac{(t_k+t_{k+1}-t_{k-1})(t_k+t_{k-1}-t_{k+1})}{4t_{k+1}t_{k-1}}}.\label{halfcoslaw}
\end{align}
The law of sines gives for $k\in\{1,2,3\}$
\begin{align}
\sin\varphi_k&=\frac{\sqrt{(t_{k+1}+t_{k-1}+t_k)(t_{k+1}+t_{k-1}-t_k)(t_k+t_{k+1}-t_{k-1})(t_k+t_{k-1}-t_{k+1})}}{2t_{k+1}t_{k-1}},\label{sinlaw}\\
\sin\left(\frac{\varphi_k}{2}\right)&=\sqrt{\frac{(t_{k+1}+t_{k-1}+t_k)(t_{k+1}+t_{k-1}-t_k)}{4t_{k+1}t_{k-1}}}.\label{halfsinlaw}
\end{align}
For $k\in\{1,2,3\}$, denoting $\dsp y_k=\cot\left(\frac{\varphi_k}{2}\right)=\frac{\dsp\cos\left(\frac{\varphi_k}{2}\right)}{\dsp\sin\left(\frac{\varphi_k}{2}\right)}=\frac{1+\cos\varphi_k}{\sin\varphi_k}=\frac{\sin\varphi_k}{1-\cos\varphi_k}$, we thus get the law of cotangents
\begin{equation}\label{cotlaw}
y_k=\sqrt{\frac{(t_k+t_{k+1}-t_{k-1})(t_k+t_{k-1}-t_{k+1})}{(t_{k+1}+t_{k-1}+t_k)(t_{k+1}+t_{k-1}-t_k)}}=\sqrt{\frac{t_k^2-(t_{k+1}-t_{k-1})^2}{(t_{k+1}+t_{k-1})^2-t_k^2}}.
\end{equation}
\begin{remark}
For $t_1=t_2=t_3$, we have $\dsp \varphi_1=\varphi_2=\varphi_3=\frac{2\pi}{3}$, so that we recover the well-known property asserting that the tension forces at the junction circle form $120^\circ$ angles.
\end{remark}
\begin{remark}
For $h=0$, the condition of criticality is given by
\[
x_1x_2x_3(t_1x_2x_3+t_2x_3x_1+t_3x_1x_2),
\]
That is
\[
x_1=0\qquad\mbox{or}\qquad x_2=0\qquad\mbox{or}\qquad x_3=0\qquad\mbox{or}\qquad\frac{t_1}{x_1}+\frac{t_2}{x_2}+\frac{t_3}{x_3}=0.
\]
As intuitively expected, we will see that only one of the first three conditions correspond to a configuration of minimum energy, namely exclusion or inclusion which happens when one of the surface tension is greater than the sum of the two others.
\end{remark}
\subsection{A more physical but equivalent way to proceed: the Lagrangian}
In fact, physics laws may appear more clearly when working (equivalently) with Lagrange multipliers as mentioned earlier, indeed, $E$ reaches its minimum in $\Rbb^4$ at a point vanishing the gradient of the Lagrangian $\mathcal{L}=E-P_1V_1-P_1V_1$ where $P_1,P_2\in\Rbb$ are Lagrange multipliers (relative pressures inside each cell). Assuming $h\neq0$ and deriving $E$ with respect to $x_1$ and $x_2$, we easily get the pressures.
\begin{align}
P_1&=-\frac{4t_1x_1}{h^2+x_1^2}=-\frac{2t_1s_1}{h}=-\frac{2t_1}{r_1},\\
P_2&=\frac{4t_2x_2}{h^2+x_2^2}=\frac{2t_2s_2}{h}=\frac{2t_2}{r_2},
\end{align}
Giving respectively the Young-Laplace law for the spherical caps $1$ and $2$ \cite{barrat,chaikin_lubensky,gennes_brochard-wyart_quere,israelachvili,safran}. Replacing $P_1$ and $P_2$ by the previous expressions in the derivatives of $\mathcal{L}$ with respect to $h$ and $x_3$, we get as previously the force balance on the junction circle \eqref{young_dupre} giving \eqref{young_dupre_angles} when expressed in terms of angles. Note that the sine relation of \eqref{young_dupre_angles} also gives the Young-Laplace law for the spherical cap $3$ and the property $t_1/r_1+t_2/r_2+t_3/r_3=0$. Finally, finding the minimum of the energy $E$ under volume constraints amounts to seek for the solutions of the following system with unknown $(x_1,x_2,x_3,h)\in\Rbb^3\times\Rbb^*$
\begin{subequations}\label{sysext}
\begin{align}
t_1\frac{h^2-x_1^2}{h^2+x_1^2}+t_2\frac{h^2-x_2^2}{h^2+x_2^2}+t_3\frac{h^2-x_3^2}{h^2+x_3^2}&=0,\label{tc}\\
t_1\frac{x_1}{h^2+x_1^2}+t_2\frac{x_2}{h^2+x_2^2}+t_3\frac{x_3}{h^2+x_3^2}&=0,\label{ts}\\
-x_1(3h^2+x_1^2)+x_3(3h^2+x_3^2)&=w_1,\label{vo}\\
x_2(3h^2+x_2^2)-x_3(3h^2+x_3^2)&=w_2,\label{vt}
\end{align}
\end{subequations}
And to seek for the minimum of the constrained energy for $h=0$. One of the points is then the global minimizer.\\
\subsection{Separation and internalization}\label{sysexindz0}
Consider a critical point $(x_1,x_2,x_3,h)$ of the constrained energy and assume that $h\neq0$. Then the force balance at the cells junction (\ref{tc},\ref{ts}) implies that the three tension forces form a triangle whose lengths are $t_1,t_2,t_3$, which can only occur if these lengths fulfill the triangle inequalities $t_k\leq t_{k+1}+t_{k-1}$ for all $k\in\{1,2,3\}$. In fact, the inequalities are strict, indeed, if $t_k=t_{k+1}+t_{k-1}$ for a $k\in\{1,2,3\}$, then the force balance implies $x_{k-1}=x_{k+1}$ which is forbidden by the volume constraints.\\
\\
Hence, if for a given $k\in\{1,2,3\}$, the tension $t_k$ satisfies $t_k\geq t_{k+1}+t_{k-1}$, then the constrained energy has no critical point in $\mcm\setminus\mcm_0$ and the global minimum, as well as any local minimum, is reached on $\mcm_0$. In fact, we can intuitively guess that the tension $t_k$ pulls until surface $k$ completely disappears, which leads to internalization of cell $1$ in cell $2$ if $k=1$, to internalization of cell $2$ in cell $1$ if $k=2$ and to separation if $k=3$. Mathematically speaking, we will prove that the constrained energy has no other minimizer than the separated/internalized configuration corresponding to the triangle inequality which is violated (mutually exclusive cases). We will also prove that these configurations are not local minimizers of the constrained energy on $\mcm$ when the triangle inequalities are fulfilled.
\begin{remark}
The three triangle inequalities $t_k<t_{k+1}+t_{k-1}$ for all $k\in\{1,2,3\}$ can be written more concisely as
\[
2\max(t_1,t_2,t_3)<t_{\mathrm{s}},
\]
Where we recall that $t_{\mathrm{s}}=t_1+t_2+t_3$. It can also be written in the following way, for a chosen $k\in\{1,2,3\}$
\[
|t_{k+1}-t_{k-1}|<t_k<t_{k+1}+t_{k-1},
\]
In this case, the similar inequalities fulfilled by the other $\ell\in\{1,2,3\}\setminus\{k\}$ are a consequence of the previous one.
\end{remark}
Keep on considering general $t_1,t_2,t_3$. First, on $\mcm_0$, the only minimizers of the energy are the points given by $x_1=0$ or $x_2=0$ or $x_3=0$, that is, the points $u_1,u_2,u_3$ defined by
\begin{subequations}
\begin{align}
u_1&=\left(0,\sqrt[3]{w_3},\sqrt[3]{w_1},0\right),\label{z10}\\
u_2&=\left(-\sqrt[3]{w_3},0,-\sqrt[3]{w_2},0\right),\label{z20}\\
u_3&=\left(-\sqrt[3]{w_1},\sqrt[3]{w_2},0,0\right).\label{z30}
\end{align}
\end{subequations}
Indeed, on $\mcm_0$, considering for instance the parameterization in $x_3$ given in \eqref{paramm0}, we have $E=\pi\psi(x_3^3)$ with $\psi:\Rbb\to\Rbb^+$ defined for $x\in\Rbb$ by
\begin{equation}\label{psi}
\psi(x)=t_1\sqrt[3]{(x-w_1)^2}+t_2\sqrt[3]{(x+w_2)^2}+t_3\sqrt[3]{x^2}.
\end{equation}
For $x\in\Rbb\setminus\{-w_2,0,w_1\}$, we have $\psi''(x)<0$ ($\psi$ is concave on each interval) and for $y\in\{-w_2,0,w_1\}$, $\lim_{x\to y^\pm}\psi'(x)=\mp\infty$. Hence, on $\mcm_0$, the energy has local minima only at the points $u_1,u_2,u_3$. Since $\dsp\lim_{x\to\pm\infty}\psi(x)=+\infty$, at least one of these points is a global minimizer of the energy on $\mcm_0$, hence also on $\mcm$ if $2\max(t_1,t_2,t_3)\geq t_{\mathrm{s}}$.\\
\\
Second, using the Hessian of the constrained energy, one can see which points are possible minimizers. Using the parameterizations of $E$ given in \eqref{paramm}, we have
\begin{align*}
\nabla_{(x_1,h)^2}E(0,0)&=2\pi\begin{pmatrix}t_1&0\\0&t_1-t_2-t_3\end{pmatrix},\\
\nabla_{(x_2,h)^2}E(0,0)&=2\pi\begin{pmatrix}t_2&0\\0&t_2-t_3-t_1\end{pmatrix},\\
\nabla_{(x_2,h)^2}E(0,0)&=2\pi\begin{pmatrix}t_3&0\\0&t_3-t_1-t_2\end{pmatrix}.
\end{align*}
Hence, for $k\in\{1,2,3\}$, $u_k$ can be a local minimizer on $\mcm$ only if $t_k\geq t_{k+1}+t_{k-1}$ and is for sure one if inequality is strict. Thus, if $2\max(t_1,t_2,t_3)<t_{\mathrm{s}}$, none of the points $u_1,u_2,u_3$ can be a local minimizer of the energy on $\mcm$, hence neither a global one. Moreover, for $2\max(t_1,t_2,t_3)\geq t_{\mathrm{s}}$, since the energy has its global minimizer among $u_1,u_2,u_3$ and since the conditions on the tensions $t_1,t_2,t_3$ are mutually exclusive, it has only one local minimizer on $\mcm$ which is also the global one.
\begin{remark}
The latter holds if for $k\in\{1,2,3\}$, $t_k=t_{k+1}+t_{k-1}$. In this case, the Hessian matrix at $u_k$ is positive but not definite, however, since the Hessian matrices at $u_{k-1}$ and $u_{k+1}$ are not positive and since the global minimizer of the energy is among $u_1,u_2,u_3$, it is necessarily $u_k$.
\end{remark}
To summarize, we have
\begin{itemize}[leftmargin=*]
\item If $2\max(t_1,t_2,t_3)<t_{\mathrm{s}}$, the energy on $\mcm$ has no local minimizer located on $\mcm_0$.

\begin{minipage}{5.7cm}\item If $t_1\geq t_2+t_3$, then $E$ reaches its global minimum at $u_1$ (which is also the unique minimizer), corresponding to the configuration of cell $1$ internalized in cell $2$. We then have $x_1=0$, $\dsp x_2=\sqrt[3]{\frac{6v_3}{\pi}}$, $\dsp x_3=\sqrt[3]{\frac{6v_1}{\pi}}$.
\end{minipage}
\hspace{0.25cm}
\begin{minipage}{5.65cm}
\begin{figure}[H]
\centering
\includegraphics[width=5.65cm]{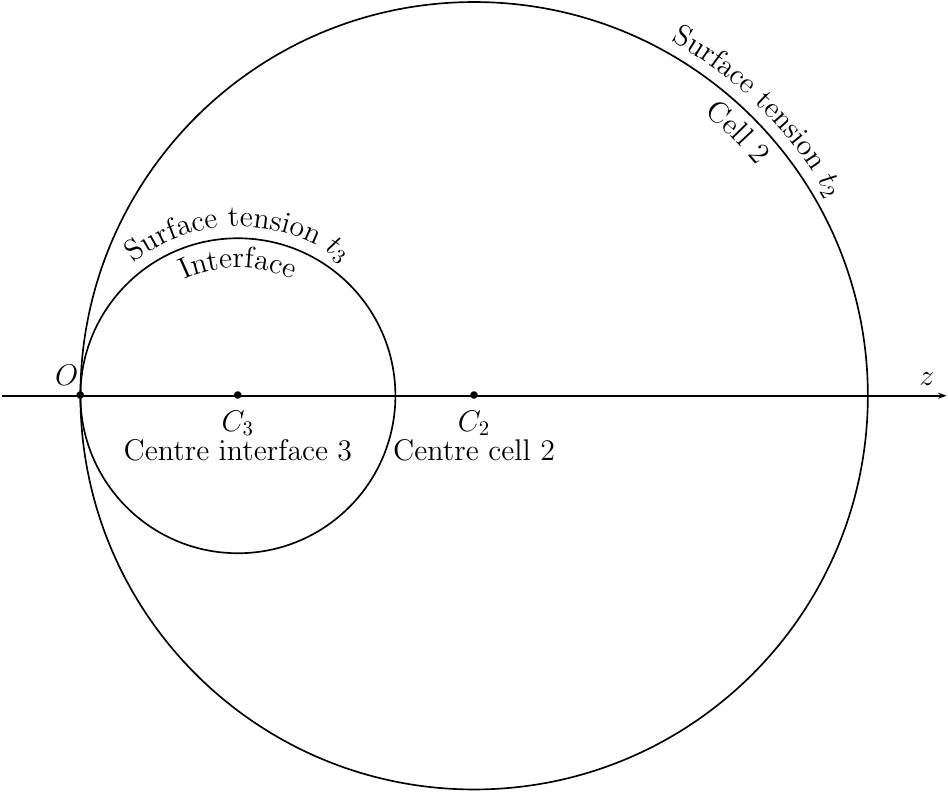}
\caption{Cell $1$ internalized in cell $2$.}\label{inter1in2}
\end{figure}
\end{minipage}

\begin{minipage}{5.7cm}\item If $t_2\geq t_1+t_3$, then $E$ reaches its global minimum at $u_2$ (which is also the unique minimizer), corresponding to the configuration of cell $2$ internalized in cell $1$. We then have $\dsp x_1=-\sqrt[3]{\frac{6v_3}{\pi}}$, $x_2=0$, $\dsp x_3=-\sqrt[3]{\frac{6v_2}{\pi}}$.
  \end{minipage}
  \hspace{0.25cm}
  \begin{minipage}{5.65cm}
  \begin{figure}[H]
    \centering
    \includegraphics[width=5.65cm]{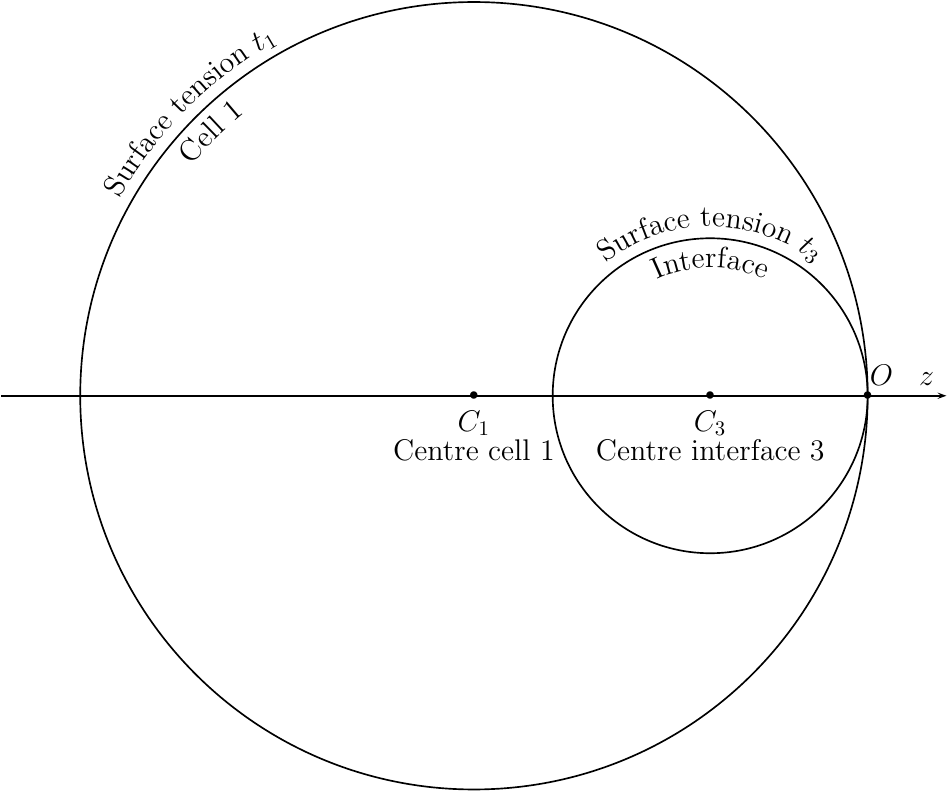}
\caption{Cell $2$ internalized in cell $1$.}\label{inter2in1}
  \end{figure}
  \end{minipage}

\begin{minipage}{5.7cm}\item If $t_3\geq t_1+t_2$, then $E$ reaches its global minimum at $u_3$ (which is also the unique minimizer), corresponding to the configuration of separated cells $1$ and $2$. We then have $\dsp x_1=-\sqrt[3]{\frac{6v_1}{\pi}}$, $\dsp x_2=\sqrt[3]{\frac{6v_2}{\pi}}$, $x_3=0$.
  \end{minipage}
  \hspace{0.25cm}
  \begin{minipage}{5.65cm}
\begin{figure}[H]
\centering
\includegraphics[width=5.65cm]{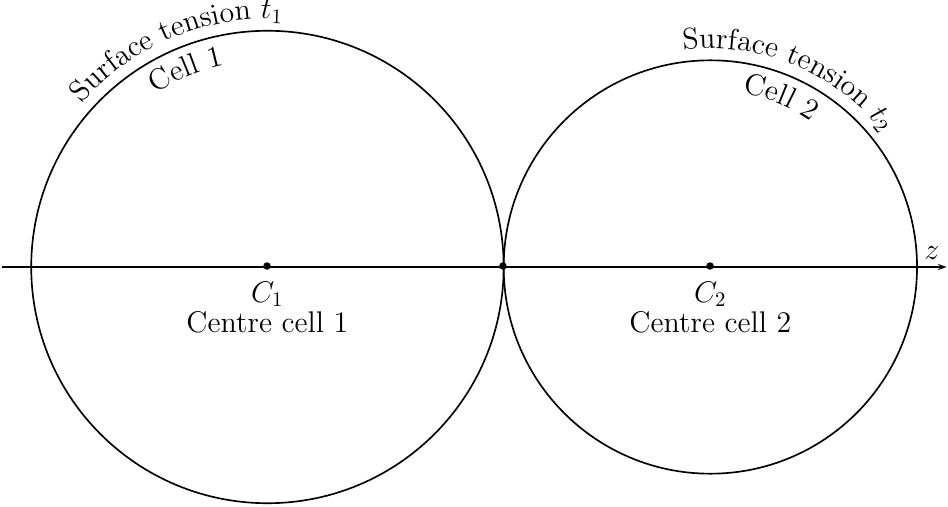}
\caption{Separated cells $1$ and $2$.}\label{sepcell}
  \end{figure}
  \end{minipage}
\end{itemize}
\begin{remark}
For all $k\in\{1,2,3\}$, if the point $x_k=0$ is not the global minimizer of the constrained energy, then it is a saddle point. The energy has two other critical points for $h=0$ that we did not take into account since they are not minimizers. As it can be seen with the function $\psi$ defined in \eqref{psi}, one of these critical point is such that $x_3\in(-\sqrt[3]{w_2},0)$ and the other one is such that $x_3\in(0,\sqrt[3]{w_1})$. In can be easily checked that these points are local maxima of the energy.
\end{remark}
\subsection{Non-degenerate case}\label{sysexindgen}
As previously mentioned, for $2\max(t_1,t_2,t_3)<t_{\mathrm{s}}$, any local minimizer of the constrained energy is a critical point in $\mcm\setminus\mcm_0$. Hence, consider a critical point of the energy in $\mcm\setminus\mcm_0$ and $k\in\{1,2,3\}$. We will prove that $|h|$ can be uniquely written in terms of $z_{k-1},z_{k+1}$, that $z_{k-1},z_{k+1}$ can be uniquely written in terms of $z_k$, and that $z_k$ is the root of a degree $5$ polynomial. We will then prove that the polynomial has a unique real root which shows that the constrained energy has a unique critical point (up to the sign of $h$ which does not change the doublet configuration) which is also the unique local and global minimizer.\\
\\
For convenience, we define a notation of the volumes, similar to the difference of angles.
\[
g_1=w_2\qquad,\qquad g_2=w_1\qquad,\qquad g_3=-w_3.
\]
First, from the volume constraints, we have
\begin{equation}\label{h}
|h|=\sqrt[3]{\frac{g_k}{z_{k+1}(3+z_{k+1}^2)-z_{k-1}(3+z_{k-1}^2)}}.
\end{equation}
Second, since $\alpha_{k\pm1}=\alpha_k\mp\varphi_{k\mp1}$, $z_{k\pm1}$ are obtained by Möbius transformations of $z_k$
\[
z_{k\pm1}=\frac{y_{k\mp1}z_k\mp1}{\pm z_k+y_{k\mp1}}.
\]
Moreover, from volume constraints, we have
\[
g_{k-1}(z_{k-1}(3+z_{k-1}^2)-z_k(3+z_k^2))=g_{k+1}(z_k(3+z_k^2)-z_{k+1}(3+z_{k+1}^2)).
\]
Hence, replacing $z_{k\pm1}$ by their expression in terms of $z_k$ we get $F_k(z_k)=0$ where the polynomial $F_k\in\Rbb[X]$ is defined by
\begin{equation}\label{popol}
F_k(X)=\frac{1+d_k}{2}T_k(X)+\frac{1-d_k}{2}U_k(X),
\end{equation}
Where $\dsp d_k=\frac{g_{k-1}-g_{k+1}}{g_{k-1}+g_{k+1}}$ and $T_k,U_k\in\Rbb[X]$ are defined by
\begin{align}\label{poltu}
T_k(X)&=\left(X+y_{k-1}\right)^3\left(\left(X-\frac{3}{2}y_{k+1}\right)^2+\frac{3}{4}y_{k+1}^2+1\right),\\
U_k(X)&=\left(X-y_{k+1}\right)^3\left(\left(X+\frac{3}{2}y_{k-1}\right)^2+\frac{3}{4}y_{k-1}^2+1\right).
\end{align}
In the following, for sake of simplicity in the notations, for a chosen $k\in\{1,2,3\}$, we simply denote $d=d_k,g=g_k,t=t_k,w=w_k,y=y_k,z=z_k,F=F_k,T=T_k,U=U_k$. The variables with subscript $k\pm1$ are denoted with a $\pm$ subscript, that is $d_{\pm}=d_{k\pm1},g_{\pm}=g_{k\pm1},t_{\pm}=t_{k\pm1},x_{\pm}=x_{k\pm1},y_{\pm}=y_{k\pm1},z_{\pm}=z_{k\pm1},F_{\pm}=F_{k\pm1},T_{\pm}=T_{k\pm1},U_{\pm}=U_{k\pm1}$.
\subsubsection{Uniqueness of the local/global minimum}\label{unique}
Since the polynomial $F$ defined in \eqref{popol} has $5$ roots in $\Cbb$, there are \textit{a priori} five critical points of the energy under volume constraints, one of them being the global minimum. However, as we are going to prove, the polynomial $F$ has a unique real root $z^*$ so that the energy has a unique critical point which is also its global minimizer. The other parameters of the configuration at minimum energy are then given by
\begin{align*}
z_\pm^*&=\frac{y_\mp z^*\mp1}{y_\mp\pm z^*},\\
\alpha^*&=2\arctan z^*,\\
\alpha_{\pm}^*&=\alpha^*\mp\varphi_{\mp},\\
h^*&=\pm\sqrt[3]{\frac{g}{z_{+}^*(3+z_{+}^{*2})-z_{-}^*(3+z_{-}^{*2})}}.
\end{align*}
\begin{theorem}\label{uniquepopol}
The polynomial $F$ defined in \eqref{popol} has a unique real root.
\end{theorem}
\begin{proof}
$T(y_\pm),U(y_\pm)\neq0$, hence for $z\in\Rbb$, $F(z)=0$ if and only if
\[
\frac{U(z)}{T(z)}=-\frac{g_{-}}{g_{+}},
\]
That is, $F(z)=0$ if and only if
\begin{equation}\label{eqfw}
f(\xi)=\frac{g_{-}}{g_{+}},
\end{equation}
Where $\dsp \xi=\frac{y_{+}-z}{y_{-}+z}$ $\left(\mbox{that is }z=\dsp\frac{y_{+}-\xi y_{-}}{1+\xi}\right)$ and the rational fraction $f$ is defined by
\begin{equation}\label{fpq}
f(\xi)=\xi^3\frac{H(\xi)}{R(\xi)},
\end{equation}
Where the polynomials $H,R\in\Rbb[X]$ (with no root in $\Rbb$) are defined by
\begin{align*}
H&=\left(y_{-}^2+1)X^2+(3y_{-}^2+y_{+}y_{-}+2\right)X+\left(y_{+}^2+3y_{-}^2+3y_{+}y_{-}+1\right),\\
R&=\left(3y_{+}^2+y_{-}^2+3y_{+}y_{-}+1\right)X^2+\left(3y_{+}^2+y_{+}y_{-}+2\right)X+\left(y_{+}^2+1\right).
\end{align*}
The derivative $f'$ of $f$ is given by
\[
f'(\xi)=\frac{S(\xi)}{R^2(\xi)}\;,\;x\in\Rbb,
\]
Where the polynomial $S\in\Rbb[X]$ is defined by
\[
S=3\delta X^2(X+1)^2\left(\left(X+\eta\right)^2+\theta\right),
\]
Where
\begin{align*}
\delta&=\left(y_{-}^2+1\right)\left(3y_{+}^2+y_{-}^2+3y_{+}y_{-}+1\right),\\
\eta&=\frac{y_{+}y_{-}\left(y_{+}^2+y_{-}^2+3y_{+}y_{-}\right)+y_{+}^2+y_{-}^2+1}{\delta},\\
\theta&=\frac{\left(y_{-}+y_{+}\right)^2\left(y_{+}^2+y_{-}^2+y_{+}y_{-}+1\right)\left(2y_{+}^2y_{-}^2+3y_{+}^2+3y_{-}^2+y_{+}y_{-}+3\right)}{\delta^2}.
\end{align*}
Since $y_{\pm}>0$, we have $\delta,\eta,\theta>0$, thus for all $\xi\in\Rbb\setminus\{-1,0\}$, we have $S(\xi)>0$ so that $f$ is a strictly increasing function. Moreover $\dsp\lim_{\xi\to\pm\infty}f(\xi)=\pm\infty$. Thus, $f$ is a bijection from $\Rbb$ to $\Rbb$, that is, for all $\zeta\in\Rbb$, there exists a unique $\xi\in\Rbb$ such that $f(\xi)=\zeta$. In particular, taking $\zeta=g_{-}/g_{+}$ proves the theorem.
\end{proof}
\begin{remark}
Theorem \ref{uniquepopol} also proves that for strictly positive volumes $v_1,v_2$ (what we assume), the unique real root of $F$ has multiplicity $1$. Indeed, the only possible roots having multiplicity larger than $1$ are $\{-1,0\}$. Since $y_++y_->0$, $\dsp z\mapsto\frac{y_{+}-z}{y_{-}+z}$ is a proper Möbius transform (not constant) and the unique respective preimages of $-1$ and $0$ are $\infty$ and $y_+$ which are not roots of $U(z)/T(z)=-g_-/g_+$.
\end{remark}
We can summarize the different cases which can appear in a double cell with constrained volumes by the fact that there always exists a unique local minimum of the energy, which is then the global minimum of the system. When the tensions fulfill the triangle inequality, the configuration of the system can be found by solving a degree $5$ equation which has a unique real root. For the cases where the tensions do not fulfill the triangle inequality, one of the surfaces disappear depending on which tension is larger than the sum of the two others. These results can be visualized on the following graph.
\begin{figure}[H]
\centering
\includegraphics[width=5cm]{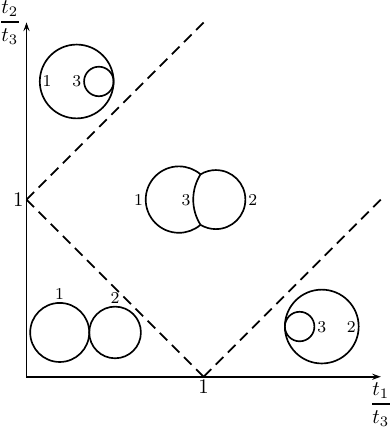}
\caption{Type of configurations depending on the tensions}\label{tensions_configurations}
\end{figure}
\begin{remark}
We cannot expect a closed formula in terms of radicals giving the roots of $F$ defined in \eqref{popol} (in particular the real one). Indeed, it is well-known that there is no solution in radicals to general polynomial equations of degree $5$ \cite{cox}. To check that $F$ is not a special polynomial which can be solved by radicals, consider $d,y_{-},y_+$ as indeterminates and $F$ as an element of $\Qbb(d,y_{-},y_+)[X]$. $T$ and $U$ are coprime in $\Cbb[d,y_{-},y_+,X]$, hence also $T-U$ and $T+U$ are. Moreover, $F=(T-U)d/2+(T+U)/2$ hence $F$ is irreducible in $\Cbb[d,y_{-},y_+,X]$. Since $F$ is monic in $X$, from Gauss lemma, $F$ is irreducible in $\Cbb(d,y_{-},y_+)[X]$, hence in $\Qbb(d,y_-,y_+)[X]$. Thus the Galois $G$ group of $F$ is a transitive subgroup of $S_5$. $G$ can be computed with a computer algebra system and we have exactly $G=S_5$ so that $F$ is not solvable by radicals. Hence, from Hilbert irreducibility theorem, polynomials obtained by specializing $d,y_{-},y_+$ in $\Qbb$ for $F$ are ``almost always'' irreducible. In fact, a consequence of this theorem is that the Galois group of the specialized polynomial is ``almost always'' isomorphic to the Galois group of $F$ over $\Qbb(d,y_{-},y_+)$ which is $S_5$ so that the specialization is not solvable by radicals.\\
Consider the simple case where $t_1=t_2=t_3$, in this case we specialize the polynomial $F$ in $y_-,y_+=1/\sqrt{3}\notin\Qbb$. Then denote $\tilde F=9\sqrt{3}F(X/\sqrt{3})$, that is
\[
\tilde F=X^5+10dX^2+15X+6d=2(5X^2+3)d+X(X^4+15).
\]
Reasoning as previously shows that $\tilde F$ is irreducible in $\Qbb(d)[X]$ and that its Galois group is $S_5$ as well, so that the polynomials obtained by specializing $d$ in $\Qbb$ for $\tilde F$ are ``almost always'' irreducible and that the Galois group of the specialized polynomial is ``almost always'' isomorphic to the Galois group of $\tilde F$ over $\Qbb(d)$ which is $S_5$, so that the specialization is not solvable by radicals. In fact, in this particular case, all specializations $d$ such that $\tilde F(d,X)$ is irreducible in $\Qbb[X]$ produce $S_5$ as Galois group. Moreover, the values of $d\in\Qbb$ for which $\tilde F(d,X)$ is reducible can be explicitly determined, they are the values of $d\in\Qbb$ for which $\tilde F(d,X)$ has a rational root. Indeed, if for $d\in\Qbb$, $\tilde F(d,X)$ factorizes as polynomials of respective degree $1$ and $4$, $\tilde F(d,X)$ has obviously a rational root. Now , if for $d\in\Qbb$, $\tilde F(d,X)$ factorizes as polynomials of respective degree $2$ and $3$ in $\Qbb[X]$, one can check that we have necessarily $d=\pm1$ (it corresponds to $w_-=0$ or $w_+=0$) and
\[
\tilde F=(X\pm1)^3(X^2\mp3X+6).
\]
In particular, $\tilde F(\pm1,X)$ has a rational root.\\
\\
The values of $d\in\Qbb$ such that $\tilde F(d,X)$ has a root in $\Qbb$ can be easily expressed since $F(d,X)$ is linear in $d$, they are
\[
d=\frac{m(m^4+15n^4)}{2n^3(5m^2+3n^2)}\qquad,\qquad m,n\in\Zbb,
\]
And in this case
\[
\tilde F=\left(X+\frac{m}{n}\right)\left(X^4-\frac{m}{n}X^3+\frac{m^2}{n^2}X^2-\frac{3m(m^2-25n^2)}{n(5m^2+3n^2)}X+\frac{3m(m^4+15n^4)}{n^2(5m^2+3n^2)}\right).
\]
It is in fact possible to get closed formulae for the roots of a general polynomial of degree $5$, one has however to use a Tschirnhaus transformation \cite{tschirnhaus_green} to write the quintic in its principal form or in its Bring-Jerrard form \cite{brioschi,hermite}. In the first case, one has to solve the associated icosahedral equation \cite{klein_morrice,nash} whose roots can be expressed in terms of hypergeometric functions \cite{abramowitz_stegun}. In the second case, the roots can be written in terms of Jacobi theta functions or hypergeometric functions \cite{abramowitz_stegun}. However, although theoretically interesting, these formulae are quite cumbersome and are often not the most appropriate way to solve the problem for numerical applications.
\end{remark}
\subsubsection{A few words about inferring the tensions}\label{teninf}
Consider the problem of expressing $t_1,t_2,t_3$ (up to a common scaling factor) in terms of the angles $\varphi_1,\varphi_2,\varphi_3$. We present here some laws which can be useful for force inference. Generally well-known, some of them are quite obvious, some others a bit more complicated. If we consider the law of sines \eqref{sinlaw}, then defining by $r$ the circumradius of the triangle defined by the tension forces, we have
\[
r=\frac{t_1t_2t_3}{\sqrt{(t_1+t_2+t_3)(t_1+t_2-t_3)(t_2+t_3-t_1)(t_3+t_1-t_2)}},
\]
So that for $k\in\{1,2,3\}$
\begin{equation}\label{lamifull}
t_k=2r\sin\varphi_k,
\end{equation}
Which is in fact Lami's theorem when expressed without the factor
\begin{equation}\label{lami}
\frac{t_1}{\sin\varphi_1}=\frac{t_2}{\sin\varphi_2}=\frac{t_3}{\sin\varphi_3}.
\end{equation}
We immediately get a similar relations for the cosines simply by writing for $k\in\{1,2,3\}$, $\sin\varphi_k=\sqrt{1-\cos^2\varphi_k}$. We can also write the tensions using the three angles and the perimeter $t_{\mathrm{s}}=t_1+t_2+t_3$ by using the following expression linking the circumradius, the perimeter and the three angles
\[
\frac{8r}{t_{\mathrm{s}}}=\frac{(\sin\varphi_1+\sin\varphi_2-\sin\varphi_3)(\sin\varphi_2+\sin\varphi_3-\sin\varphi_1)(\sin\varphi_3+\sin\varphi_1-\sin\varphi_2)}{\sin^2\varphi_1\sin^2\varphi_2\sin^2\varphi_3}
\]
We then get for $k\in\{1,2,3\}$
\begin{equation}
t_k=\frac{t_{\mathrm{s}}\beta}{4\sin\varphi_{k+1}\sin\varphi_{k-1}},
\end{equation}
Where
\[
\beta=\frac{(\sin\varphi_1+\sin\varphi_2-\sin\varphi_3)(\sin\varphi_2+\sin\varphi_3-\sin\varphi_1)(\sin\varphi_3+\sin\varphi_1-\sin\varphi_2)}{\sin\varphi_1\sin\varphi_2\sin\varphi_3}.
\]
Lami's theorem again follows by dividing each $t_k$ by $\sin\varphi_k$ for $k\in\{1,2,3\}$. A similar expression for the cosines immediately follows, but we can also get it from the law of cosines \eqref{coslaw} and we have for $k\in\{1,2,3\}$
\[
t_k=\frac{t_{\mathrm{s}}}{2}\left(\frac{1}{1-\cos\varphi_{k+1}}+\frac{1}{1-\cos\varphi_{k-1}}-\frac{1-\cos\varphi_k}{(1-\cos\varphi_{k+1})(1-\cos\varphi_{k-1})}\right).
\]
This gives in terms of sines of the half angles, for $k\in\{1,2,3\}$
\[
t_k=\frac{t_{\mathrm{s}}}{4}\left(\frac{1}{\dsp\sin^2\left(\frac{\varphi_{k+1}}{2}\right)}+\frac{1}{\dsp\sin^2\left(\frac{\varphi_{k-1}}{2}\right)}-\frac{\dsp\sin^2\left(\frac{\varphi_k}{2}\right)}{\dsp\sin^2\left(\frac{\varphi_{k+1}}{2}\right)\sin^2\left(\frac{\varphi_{k-1}}{2}\right)}\right).
\]
The relation for the cosines of the half angles immediately follows. We can also write the tensions in terms of the perimeter using the law of cotangents \eqref{cotlaw} and we have all $k\in\{1,2,3\}$.
\[
t_k=\frac{t_{\mathrm{s}}}{2}\cot\left(\frac{\varphi_k}{2}\right)\left[\cot\left(\frac{\varphi_{k+1}}{2}\right)+\cot\left(\frac{\varphi_{k-1}}{2}\right)\right].
\]
If we now consider the problem of expressing the tensions in terms of the interfaces radii and the distance between their centers, taking the convention $h\in\Rbb^+$ with signed radii $r_k,k\in\{1,2,3\}$ as defined in \eqref{zar}, from equation \eqref{lamifull}, we get for all $k\in\{1,2,3\}$
\[
t_k=\frac{hr(C_{k+1}-C_{k-1})}{r_{k+1}r_{k-1}},
\]
That is, with unsigned radii $r_k,k\in\{1,2,3\}$
\[
t_k=\frac{hr\left|C_{k+1}-C_{k-1}\right|}{r_{k+1}r_{k-1}},
\]
So that
\begin{equation}
\frac{t_1}{r_1\left|C_2-C_3\right|}=\frac{t_2}{r_2\left|C_3-C_1\right|}=\frac{t_3}{r_3\left|C_3-C_1\right|}
\end{equation}
\subsubsection{Configuration of the cell doublet for given pressures}
We propose now to express the parameters of a cell doublet when the surface tensions $t_1,t_2,t_3$ verifying the triangle equation and when the pressures $P_1,P_2>0$ are known. Then from the Young-Laplace law and the sine relation of the Young-Dupré law
\begin{align}\label{sinerel}
  \frac{4t_1x_1}{h^2+x_1^2}=-P_1,\\
  \frac{4t_2x_2}{h^2+x_2^2}=P_2,\\
  \frac{4t_3x_3}{h^2+x_3^2}=P_3.\\
\end{align}
Where $P_3=P_1-P_2$. Hence, the cosine relation of the Young-Dupré law implies that
\begin{equation}
\pm\sqrt{4t_1^2-P_1^2h^2}\pm\sqrt{4t_2^2-P_2^2h^2}\pm\sqrt{4t_3^2-P_3^2h^2}=0.
\end{equation}
Hence
\begin{equation}
h^2=\frac{(t_1+t_2+t_3)(t_2+t_3-t_1)(t_3+t_1-t_2)(t_1+t_2-t_3)}{(P_1t_2-P_2t_1)^2+P_1P_2(t_1-t_2+t_3)(t_2-t_1+t_3)}.
\end{equation}
Since $t_1,t_2,t_3$ fulfill the triangle inequality, the right hand side is always positive. Choosing $h>0$, we then have
\begin{equation}
h=\sqrt{\frac{(t_1+t_2+t_3)(t_2+t_3-t_1)(t_3+t_1-t_2)(t_1+t_2-t_3)}{(P_1t_2-P_2t_1)^2+P_1P_2(t_1-t_2+t_3)(t_2-t_1+t_3)}}.
\end{equation}
The $x$ parameters can then be obtained by solving the sine relations \eqref{sinerel}, giving
\begin{align}
x_1&=\frac{-2t_1\pm \sqrt{4t_1^2-P_1^2h^2}}{P_1},\\
x_2&=\frac{2t_2\pm \sqrt{4t_2^2-P_2^2h^2}}{P_2},\\
x_3&=\frac{2t_3\pm \sqrt{4t_3^2-P_3^2h^2}}{P_3}.
\end{align}
Note that we do not have continuity to $0$ when $P_1=P_2$ in the coordinate $x_3$ if the minus sign is chosen, we have instead $x_3\to\pm\infty$. Remark also that the signs can not be chosen independently since $x_1,x_2,x_3$ have to fulfill the cosine relation
\begin{equation}
t_1\frac{h^2-x_1^2}{h^2+x_1^2}+t_2\frac{h^2-x_2^2}{h^2+x_2^2}+t_3\frac{h^2-x_3^2}{h^2+x_3^2}=0.
\end{equation}
More precisely, denoting by
\begin{equation}
\Delta=\sqrt{(P_1t_2-P_2t_1)^2+P_1P_2(t_1-t_2+t_3)(t_2-t_1+t_3)},
\end{equation}
We have two possible triples of solutions for $x$, denoted by $(x_1^-,x_2^-,x_3^-)$ and $(x_1^+,x_2^+,x_3^+)$. They are given by
\begin{align}
x_1^\pm&=-\frac{2t_1\pm\frac{\dsp P_1(t_1^2+t_2^2-t_3^2)-2P_2t_1^2}{\dsp\Delta}}{P_1},\\
x_2^\pm&=\frac{\dsp 2t_2\pm\frac{P_2(t_1^2+t_2^2-t_3^2)-2P_1t_2^2}{\dsp \Delta}}{P_2},\\
x_3^\pm=&\frac{2t_3\pm\frac{\dsp P_1(t_2^2+t_3^2-t_1^2)+P_2(t_1^2+t_3^2-t_2^2)}{\dsp \Delta}}{P_1-P_2}.
\end{align}
Where the possible solution $x^+$ is the one such that the third coordinate lacks continuity to $0$ as $P_1=P_2$. We will see that the triple $(x_1^+,x_2^+,x_3^+)$ does not solve the problem anyway. To solve it, the first coordinate of the triple should be negative, the second one should be positive and the third one between both. We already know that $x_1^\pm<0$ and $x_2^\pm>0$ since they respectively have the sign of $-P_1$ and $P_2$. Now, we have
\begin{align}
  x_1^--x_3^-&=\frac{t_1(P_1-P_2)+t_3P_1-\Delta}{P_1(P_1-P_2)},\\
  x_2^--x_3^-&=\frac{t_2(P_2-P_1)+t_3P_2-\Delta}{P_2(P_2-P_1)},\\
  x_1^+-x_3^+&=\frac{t_1(P_1-P_2)+t_3P_1+\Delta}{P_1(P_1-P_2)},\\
  x_2^+-x_3^+&=\frac{t_2(P_2-P_1)+t_3P_2+\Delta}{P_2(P_2-P_1)}.
\end{align}
Moreover
\begin{align}
(t_1(P_1-P_2)+t_3P_1)^2-\Delta^2=P_1(t_1+t_2+t_3)(t_1-t_2+t_3)(P_1-P_2),\\
(t_2(P_2-P_1)+t_3P_2)^2-\Delta^2=P_2(t_1+t_2+t_3)(t_2-t_1+t_3)(P_2-P_1).\\
\end{align}
Hence, we always have $x_1^-<x_3^-<x_2^-$ whereas when $P_1>P_2$, we have $x_3^+>x_2^+$ and when $P_2>P_1$, we have $x_1^+>x_3^+$. Thus, for all positive surface tensions $t_1,t_2,t_3$ fulfilling the triangle inequality and for all positive pressures $P_1,P_2$, there exists a configuration of the cell doublet, this configuration is unique and given by the parameters $x_1^-,x_2^-,x_3^-,h$ previously computed.
\section{Cell doublet with line tension}\label{cellconfminenlineten}
We consider the same problem as in section \ref{cellconfminen} except that we now assume that the system has a line tension along the junction circle. Then, the energy is now defined as
\[
E=\pi(t_1x_1^2+t_2x_2^2+t_3x_3^2+t_{\mathrm{s}} h^2+2\kappa h),
\]
Where $\kappa$ is the line tension and we recall, $t_{\mathrm{s}}=t_1+t_2+t_3$. We now restrict the problem to $h\geq 0$ (otherwise we can put an absolute value in the energy). $E$ is positive and $\dsp\lim_{(x_1,x_2,x_3,h)\to+\infty}E(x_1,x_2,x_3,h)=+\infty$, moreover, the semialgebraic set $\mcm_+=\{(x_1,x_2,x_3,h)\in\mcm\;:\;h\geq0\}$ is closed, hence $E$ reaches its minimum on $\mcm^+$, either at a critical point on $\mcm_+^*=\{(x_1,x_2,x_3,h)\in\mcm\;:\;h>0\}$ or on $\mcm_0$.
\subsection{Equations for the critical points of the constrained energy}\label{sysexindlineten}
With the same notations as in section \ref{sysexind}, the critical ideal $I$ of the energy on the affine variety $\mcm$ is given by
\begin{equation}
  \begin{aligned}
I=&I_{\mcm}+\left\langle t_1X_1(Y^2+X_2^2)(Y^2+X_3^2)+t_2X_2(Y^2+X_3^2)(Y^2+X_1^2)+\right.\\
&\left.t_3X_3(Y^2+X_1^2)(Y^2+X_2^2),2t_2X_2Y(X_1-X_2)(Y^2+X_3^2)+\right.\\
&\left.2t_3X_3Y(X_1-X_3)(Y^2+X_2^2)+(t_{\mathrm{s}}Y+\kappa)(Y^2+X_2^2)(Y^2+X_3^2),\right.\\
&\left.2t_3X_3Y(X_2-X_3)(Y^2+X_1^2)+2t_1X_1Y(X_2-X_1)(Y^2+X_3^2)+\right.\\
&\left.(t_{\mathrm{s}}Y+\kappa)(Y^2+X_3^2)(Y^2+X_1^2),2t_1X_1Y(X_3-X_1)(Y^2+X_2^2)+\right.\\
&\left.2t_2X_2Y(X_3-X_2)(Y^2+X_1^2)+(t_{\mathrm{s}}Y+\kappa)(Y^2+X_1^2)(Y^2+X_2^2)\right\rangle.
\end{aligned}
\end{equation}
The algebraic degree of the problem is $83$. We can see that the energy has no critical point such that $h=0$. However, since $\mcm_0$ is now the boundary of the domain of the constrained energy, points such that $h=0$ do not need to be critical points to be local extrema. The zeros $(x_1,x_2,x_3,h)$ of the critical ideal fulfill the modified Young-Dupré law
\begin{align}\label{modified_young_dupre}
t_1\frac{h^2-x_1^2}{h^2+x_1^2}+t_2\frac{h^2-x_2^2}{h^2+x_2^2}+t_3\frac{h^2-x_3^2}{h^2+x_3^2}+\frac{\kappa}{h}&=0,\\
t_1\frac{x_1}{h^2+x_1^2}+t_2\frac{x_2}{h^2+x_2^2}+t_3\frac{x_3}{h^2+x_3^2}&=0,
\end{align}
And the volumes constraints. One can again work the Lagrangian of the energy $\mathcal{L}=E-P_1V_1-P_2V_2$ to recover this results. Making the derivatives of $\mathcal{L}$ with respect to $x_1,x_2,x_3,h$ vanish, we get the Laplace law on the surfaces $1,2,3$
\begin{align}
P_1&=-\frac{4t_1x_1}{h^2+x_1^2}=-\frac{2t_1s_1}{h}=-\frac{2t_1}{r_1},\\
P_2&=\frac{4t_2x_2}{h^2+x_2^2}=\frac{2t_2s_2}{h}=\frac{2t_2}{r_2},\\
P_1-P_2&=\frac{4t_2x_3}{h^2+x_3^2}=\frac{2t_3s_3}{h}=\frac{2t_3}{r_3},
\end{align}
And the modified Young-Dupré law previously mentioned. Before tackling the local extrema of the energy such that $h\neq0$, consider first the case of points such that $h=0$ (points in $\mcm_0$). Points which are local minima of the energy on $\mcm_0$ are of course the same as in the problem \ref{sysexindz0} without line tension and they are those corresponding to internalization or separation. In the presence of line tension, the three points are local minima of the energy, indeed, using the parameterizations of $E$ given in \eqref{paramm}, we have
\begin{align}
E(x_1,h)&\mathop{=}_{x_1,h\to0}\pi\left(t_3w_1^{2/3}+t_2w_3^{2/3}+t_1x_1^2+2\kappa h\right)+o(t_1x_1^2+2\kappa h),\\
E(x_2,h)&\mathop{=}_{x_2,h\to0}\pi\left(t_3w_2^{2/3}+t_1w_3^{2/3}+t_2x_2^2+2\kappa h\right)+o(t_2x_2^2+2\kappa h),\\
E(x_3,h)&\mathop{=}_{x_3,h\to0}\pi\left(t_1w_1^{2/3}+t_2w_2^{2/3}+t_3x_3^2+2\kappa h\right)+o(t_3x_3^2+2\kappa h).
\end{align}
We now want to find the local extrema of the energy of for $h>0$. We need to solve the following system with known $(x_1,x_2,x_3,h)$ to find the critical points of the constrained energy
\begin{align}
t_1\frac{h^2-x_1^2}{h^2+x_1^2}+t_2\frac{h^2-x_2^2}{h^2+x_2^2}+t_3\frac{h^2-x_3^2}{h^2+x_3^2}+\frac{\kappa}{h}&=0,\\
t_1\frac{x_1}{h^2+x_1^2}+t_2\frac{x_2}{h^2+x_2^2}+t_3\frac{x_3}{h^2+x_3^2}&=0,\\
-x_1(x_1^2+3h^2)+x_3(x_3^2+3h^2)&=w_1,\\
x_2(x_2^2+3h^2)-x_3(x_3^2+3h^2)&=w_2.
\end{align}
In the case without line tension, the system is easier to solve because the three forces form a triangle whose angles are uniquely defined by its sides (surface tensions). However, we now have four forces forming a quadrilateral whose angles are not uniquely defined by its sides so that the two constraints equations play a role in the angles. As previously explained, we want to be sure to find the global minimum of the energy and the local minima as well, so the system can not be solved with Newton method. Computer Algebra Systems show some difficulties in finding a Gröbner basis for this system but as previously mentioned, they are often not the best solution. Hence, instead, we propose to solve the system using Polynomial Homotopy Continuation and to this end we use \href{http://homepages.math.uic.edu/~jan/download.html}{PHCpack} \cite{verschelde}. We give the following polynomial system to the software
\begin{subequations}\label{linesys}
\begin{align}
P(z_1,z_2,z_3)\left(t_1\frac{1-z_1^2}{1+z_1^2}+t_2\frac{1-z_2^2}{1+z_2^2}+t_3\frac{1-z_3^2}{1+z_3^2}+\kappa y\right)&=0,\label{linesys1}\\
P(z_1,z_2,z_3)\left(t_1\frac{z_1}{1+z_1^2}+t_2\frac{z_2}{1+z_2^2}+t_3\frac{z_3}{1+z_3^2}\right)&=0\label{linesys2},\\
-z_1(z_1^2+3)+z_3(z_3^2+3)-w_1y^3&=0,\label{linesys3}\\
z_2(z_2^2+3)-z_3(z_3^2+3)-w_2y^3&=0,\label{linesys4}
\end{align}
\end{subequations}
With unknowns $z_1,z_2,z_3$, $y=1/h$ and where
\begin{equation}\label{polz}
P=(1+z_1^2)(1+z_2^2)(1+z_3^2).
\end{equation}
 One can also define
\begin{equation}\label{redvar}
w=(w_2-w_1)/w_3\quad,\quad\rho=y\sqrt[3]{w_3}\quad,\quad\tau_k=t_k\sqrt[3]{w_3}/\kappa\;\;\mbox{for}\;\;k\in\{1,2,3\},
\end{equation}
And solve the following system of unknowns $z_1,z_2,z_2,\rho$ instead
\begin{equation}\label{redlinesys}
\begin{aligned}
P(z_1,z_2,z_3)\left(\tau_1\frac{1-z_1^2}{1+z_1^2}+\tau_2\frac{1-z_2^2}{1+z_2^2}+\tau_3\frac{1-z_3^2}{1+z_3^2}+\rho\right)&=0,\\
P(z_1,z_2,z_3)\left(\tau_1\frac{z_1}{1+z_1^2}+\tau_2\frac{z_2}{1+z_2^2}+\tau_3\frac{z_3}{1+z_3^2}\right)&=0,\\
-2z_1(z_1^2+3)+2z_3(z_3^2+3)-(1-w)\rho^3&=0,\\
2z_2(z_2^2+3)-2z_3(z_3^2+3)-(1+w)\rho^3&=0.
\end{aligned}
\end{equation}
\begin{remark}\label{eqw}
$\rho$ can be easily removed from the unknowns by using the value given by the first equation and replace $\rho$ in the two last equations. One then works with the three last equations and the unknowns $z_1,z_1,z_2$. Note also that the two last equations can also be written as
\begin{align*}
z_2(z_2^2+3)-z_1(z_1^2+3)-\rho^3&=0,\\
(1+w)z_1(z_1^2+3)+(1-w)z_2(z_2^2+3)-2z_3(z_3^2+3)&=0.
\end{align*}
\end{remark}
Once the critical points are known, we need to determine if they are local minima, saddle points or local maxima. To this end, we use the Hessian of the Lagrangian $\mathcal{L}$ with respect to $(x_1,x_2,x_3,h)$. Indeed, consider a critical point $(x_1,x_2,x_3,h)$, its associated Lagrange multipliers $P_1,P_2$ and $H$ the Hessian of $\mathcal{L}$ with respect to the geometry taken at $(x_1,x_2,x_3,h,P_1,P_2)$. Denote $T$ the tangent space to $\mcm$ at $(x_1,x_2,x_3,h)$, that is $T=\langle\nabla V_1,\nabla V_2\rangle^*$, where $*$ denotes the orthogonal.
\begin{itemize}
\item If $(x_1,x_2,x_3,h)$ is a local maximum of the constrained energy, then for all $d\in T,d^*Hd\leq0$.\\
If for all $d\in T\setminus\{0\},d^*Hd<0$, then $(x_1,x_2,x_3,h)$ is a local maximum of the constrained energy.
\item If $d^*Hd$ changes sign for $d\in T$, then $(x_1,x_2,x_3,h)$ is a saddle point of the constrained energy.
\item If $(x_1,x_2,x_3,h)$ is a local minimum of the constrained energy, then for all $d\in T,d^*Hd\geq0$\\
If for all $d\in T\setminus\{0\},d^*Hd>0$, then $(x_1,x_2,x_3,h)$ is a local minimum of the constrained energy.
\item If the Hessian is negative but not definite on the tangent space $T$, more analysis is required to know if $(x_1,x_2,x_3,h)$ is a local maximum or a saddle point. Similarly, if the Hessian is positive but not definite on the tangent space $T$, more analysis is required to know if $(x_1,x_2,x_3,h)$ is a local minimum or a saddle point.
\end{itemize}
Comparing the energy at the local minima and at the points $h=0,x_k=0$, $k\in\{1,2,3\}$, one can find the global minimum of the constrained energy. To test the extremality of the critical points, we need the tangent space $T$ and the Hessian matrix $H$. We recall that
\begin{equation}
\nabla V_1=\frac{\pi}{2}\begin{pmatrix}-(h^2+x_1^2)\\0\\h^2+x_3^2\\2h(x_3-x_1)\end{pmatrix}\qquad,\qquad\nabla V_2=\frac{\pi}{2}\begin{pmatrix}0\\h^2+x_2^2\\-(h^2+x_3^2)\\2h(x_2-x_3)\end{pmatrix}.
\end{equation}
Hence the following vectors form a basis of $T$.
\begin{equation}
U=\begin{pmatrix}\dsp1+c_1\\\dsp1+c_2\\\dsp1+c_3\\0\end{pmatrix}\qquad,\qquad W=\begin{pmatrix}\dsp-s_1\\\dsp-s_2\\\dsp-s_3\\1\end{pmatrix}.
\end{equation}
The Hessian matrix $H$ is simply given by
\begin{equation}
H=2\pi\begin{pmatrix}
\dsp t_1c_1&0&0&\dsp-t_1s_1\\
0&\dsp t_2c_2&0&\dsp-t_2s_2\\
0&0&\dsp t_3c_3&\dsp-t_3s_3\\
\dsp-t_1s_1&\dsp-t_2s_2&\dsp-t_3s_3&\dsp t_1c_1+t_2c_2+t_3c_3
\end{pmatrix}.
\end{equation}
Hence, the Hessian matrix $H_T$ in the basis $(U,W)$ writes
\begin{equation}\label{hessht}
H_T=2\pi\left[M_H-\kappa y\begin{pmatrix}
\dsp1&\dsp0\\
\dsp0&\dsp1
\end{pmatrix}\right],
\end{equation}
Where
\begin{equation}
M_H=\begin{pmatrix}
\dsp\sum_{k=1}^3t_kc_k^2(2+c_k)&\dsp-\sum_{k=1}^3t_kc_ks_k(2+c_k)\\
\dsp\dsp\dsp-\sum_{k=1}^3t_kc_ks_k(2+c_k)&\dsp\sum_{k=1}^3t_ks_k^2(2+c_k)
\end{pmatrix}.
\end{equation}
Thus, the $2\times2$ matrix $H_T$ is positive if and only if its trace and its determinant are positive. The trace $\Tr(H_T)$ of $H_T$ is given by
\begin{equation}\label{trhessht}
\Tr(H_T)=2\pi(2t_{\mathrm{s}}-3\kappa y),
\end{equation}
And its determinant $\det(H_T)$ by
\begin{equation}\label{dethessht}
\det(H_T)=4\pi^2(\det(M_H)-2\kappa y(t_{\mathrm{s}}-\kappa y)),
\end{equation}
Where the determinant $\det(M_H)$ of $M_H$ is given by
\begin{equation}
\det(M_H)=\left(\sum_{k=1}^3t_kc_k^2(2+c_k)\right)\left(\sum_{k=1}^3t_ks_k^2(2+c_k)\right)-\left(\sum_{k=1}^3t_kc_ks_k(2+c_k)\right)^2.
\end{equation}
We recall that the three local minima $u_1,u_2,u_3$ of $E$ for $h=0$ are given by \eqref{z10}, \eqref{z20}, \eqref{z30}
\begin{align}
u_1&=\left(0,\sqrt[3]{w_3},\sqrt[3]{w_1},0\right),\\
u_2&=\left(-\sqrt[3]{w_3},0,-\sqrt[3]{w_2},0\right),\\
u_3&=\left(-\sqrt[3]{w_1},\sqrt[3]{w_2},0,0\right),
\end{align}
The energy at these points is
\begin{align}
E(u_1)&=\pi\left(t_2\sqrt[3]{w_3^2}+t_3\sqrt[3]{w_1^2}\right),\\
E(u_2)&=\pi\left(t_1\sqrt[3]{w_3^2}+t_3\sqrt[3]{w_2^2}\right),\\
E(u_3)&=\pi\left(t_1\sqrt[3]{w_2^2}+t_2\sqrt[3]{w_2^2}\right).
\end{align}
\subsection{Theoretical and numerical results, observations}
Remark first that like in the case of surface tensions only, different relations can be obtained from the balances of forces. More precisely, we have the following lemma
\begin{lemma}
Assume that a volume constrained cell doublet with tensions $t_1,t_2,t_3,\kappa$ and geometric parameters $(z_1,z_2,z_3,y)$ is at a local minimum of energy, then we recall that equations \eqref{linesys1} and \eqref{linesys2} are fulfilled (balances of forces).
\begin{align*}
t_1c_1+t_2c_2+t_3c_3+\kappa y&=0,\\
t_1s_1+t_2s_2+t_3s_3&=0.
\end{align*}
We recall that for $k\in\{1,2,3\}$, we denote $\varphi_k=\alpha_{k+1}-\alpha_{k-1}$. Then for all $k\in\{1,2,3\}$
\begin{subequations}
\begin{align}
t_k+t_{k+1}\cos\varphi_{k-1}+t_{k-1}\cos\varphi_{k+1}+\kappa y c_k&=0,\label{cosrelk}\\
t_{k+1}\sin\varphi_{k-1}-t_{k-1}\sin\varphi_{k+1}+\kappa y s_k&=0,\label{sinrelk}.
\end{align}
\end{subequations}
In particular, from \eqref{cosrelk}, we get
\begin{equation}
t_1^2+t_2^2+t_3^2+2t_2t_3\cos\varphi_1+2t_3t_1\cos\varphi_2+2t_1t_2\cos\varphi_3=\kappa^2y^2,
\end{equation}
And for $k\in\{1,2,3\}$
\begin{equation}
\cos\varphi_k=\frac{t_k^2-t_{k-1}^2-t_{k+1}^2+\kappa y(t_kc_k-t_{k-1}c_{k-1}-t_{k+1}c_{k+1})}{2 t_{k-1}t_{k+1}},
\end{equation}
That is
\begin{equation}
\cos\varphi_k=\frac{t_k^2-t_{k-1}^2-t_{k+1}^2+\kappa y(2t_kc_k+\kappa y)}{2t_{k-1}t_{k+1}},
\end{equation}
\end{lemma}
\begin{proof}
The relations are simply obtained by linear combination of equations \eqref{linesys1} and \eqref{linesys2} with cosine/sine coefficients.
\end{proof}
\subsubsection{Inequalities for tensions and volumes}
We first propose to show an obvious inequality between the tensions and the inverse of the junction height $y$.
\begin{lemma}
If a cell doublet with tensions $t_1,t_2,t_3,\kappa$ and geometric parameters $(z_1,z_2,z_3,y)$ is at a local minimum of energy, then, the balance of forces (equations \eqref{linesys1}, \eqref{linesys2}) implies that three of the triangle inequalities the triangle inequality of the doublet without line tension is replaced by the quadrilateral inequality.
\begin{equation}\label{ineqkappa}
2\max(t_1,t_2,t_3,\kappa y)\leq t_1+t_2+t_3+\kappa y
\end{equation}
\end{lemma}
\begin{proof}
From equation \eqref{linesys1}, we have
\[-\kappa y=t_1c_1+t_2c_2+t_3c_3.\]
Hence
\[\kappa y\leq t_1+t_2+t_3\]
For $k\in\{1,2,3\}$, the inequality
\[t_k\leq\kappa y+t_{k+1}+t_{k-1}\]
Is obtained similarly by using \eqref{cosrelk}.
\end{proof}
We now propose to show a necessary condition on the tensions and volumes in order to have solutions to the system \eqref{linesys}.
\begin{lemma}
A necessary condition on the system \eqref{linesys} to have a real solution is
\begin{equation}
u\geq3\qquad\mbox{where}\qquad u=\frac{t_{\mathrm{s}}}{\kappa}\frac{\sqrt[3]{w_3}}{\sqrt[3]{2}}.
\end{equation}
Moreover, if $(z_1,z_2,z_3,y)$ is a real solution to \eqref{linesys} with $y>0$, then
\begin{equation}\label{yineq}
y\geq\frac{\sqrt[3]{4}}{\sqrt[3]{w_3}}.
\end{equation}
If for $k\in\{1,2,3\}$ we have $t_k-t_{k+1}-t_{k-1}\geq0$, then another necessary condition on the system \eqref{linesys} to have a real solution is
\begin{equation}
t_{k+1}+t_{k-1}\geq t_{\mathrm{s}}\frac{\dsp\cosh\left(\frac{1}{3}\arcsinh u^3\right)}{\sqrt{u^6+1}}
\end{equation}
\end{lemma}
\begin{proof}
Assume that the system \eqref{linesys} has a real solution $(z_1,z_2,z_3,y)$. For $k\in\{1,2,3\}$, denote $q_k=z_k(z_k^2+3)$, then $c_k=f(q_k/2)$ where for all $x\in\Rbb$
\[
f(x)=\frac{\dsp\sqrt[3]{\sqrt{x^2+1}+x}+\sqrt[3]{\sqrt{x^2+1}-x}}{\dsp\sqrt{x^2+1}}-1=2\frac{\dsp\cosh\left(\frac{1}{3}\arcsinh x\right)}{\sqrt{x^2+1}}-1.
\]
$f$ is symmetric, increasing on $\Rbb^-$ and decreasing on $\Rbb^+$. Moreover, for all $k\in\{1,2,3\}$, $|q_k|\leq w_3y^3$, hence
\[
t_1c_1+t_2c_2+t_3c_3+\kappa y\geq t_{\mathrm{s}}f(w_3y^3/2)+\kappa y,
\]
Where we recall that $t_{\mathrm{s}}=t_1+t_2+t_3$. Hence
\[
t_{\mathrm{s}}f(w_3y^3/2)+\kappa y\leq0.
\]
$f$ is positive on $(0,2)$ so that $y$ necessarily fulfills \eqref{yineq}. Denote by $\dsp \omega=y\sqrt[3]{\frac{w_3}{2}}$ and by $g$ the function $x\mapsto -f(x^3)/x$. Since $t_1c_1+t_2c_2+t_3c_3+\kappa y=0$, then $g(\omega)\geq a$ where
\[
a=\frac{\kappa}{t_{\mathrm{s}}}\sqrt[3]{\frac{2}{w_3}}.
\]
Hence $\dsp\max_{x\in\Rbb_+^*}g\geq a$. Moreover, $g$ reaches its maximum at $\omega_0$ where
\[
\omega_0=\sqrt[6]{32+\frac{4\sqrt{778}}{\sqrt{3}}\cos\left(\frac{1}{3}\arccos\left(\frac{50115\sqrt{3}}{3112\sqrt{778}}\right)\right)}.
\]
Thus $a\leq M$, with $M=g(\omega_0)$, that is
\[
\frac{\kappa}{t_{\mathrm{s}}}\sqrt[3]{\frac{2}{w_3}}\leq M,
\]
With
\[
M=\sqrt[6]{\frac{4}{3}}\sqrt[6]{2\sqrt{5}\sqrt{10883}\cos\left(\frac{1}{3}\arccos\left(\frac{1}{8}\frac{101454517}{5^{3/2}10883^{3/2}}\right)-\frac{\pi}{3}\right)-239}
\]
That is $\dsp M\simeq0.3217\simeq\frac{1}{3}$. Finally
\[
\frac{\kappa}{t_{\mathrm{s}}}\leq\frac{\sqrt[3]{w_3}}{3\sqrt[3]{2}}.
\]
Prove now the second inequality on the tensions. Consider $k\in\{1,2,3\}$ and assume that $t_k-t_{k+1}-t_{k-1}\geq0$. Then, using \eqref{cosrelk}, we have
\[
t_k-t_{k+1}-t_{k-1}+\kappa y f(w_3y^3/2)\leq0,
\]
The function $x\mapsto xf(x^3)$ is decreasing on $[\sqrt[3]{2},+\infty)$, in addition, we have $y\leq t_{\mathrm{s}}/\kappa$, hence
\[
t_k-t_{k+1}-t_{k-1}\leq-t_{\mathrm{s}}f(u^3),
\]
So that the second inequality is proven.
\end{proof}
\subsubsection{Uniqueness of the interior minimum}\label{unilocmin}
After computing the critical points for a grid of parameters $t_1,t_2,t_3,k,w_1,w_2$, there seems to be at most six critical points in the considered domain, more precisely, there are $0$, $2$, $4$ or $6$. However, it seems that there is no local minimum when there are less than six critical points, whereas there is exactly one minimum when there are six critical points. Consequently, since when exists, the interior minimum is unique, a minimization algorithm could fit to search for the latter. However, there could be some difficulties to avoid the algorithm to converge to one of the three minima on the boundary of the domain ($h=0$). This could probably be avoided by an appropriate choice of the initial guess. One can also remark that even when the energy has an interior minimum, this one is not necessarily the global minimum of the system. In these cases, the global minimum is reached when $h=0$, that is when $x_k=0$ for one of $k\in\{1,2,3\}$.
\subsubsection{Energy and singularities}
The singularities of the system are the critical points such that the determinant of the matrix \eqref{hessht} is $0$, that is, the solutions to the system \eqref{linesys} such that Jacobian has determinant $0$. For sake of simplicity, considering the reduced system \eqref{redlinesys}, one can choose one of the parameters as being a new unknown (for instance $w$) and add the equation of vanishing Jacobian determinant to the system. Solving this new system for chosen tensions $\tau_1,\tau_2,\tau_3$ gives the singular points $(z_1,z_2,z_3,\rho)$ for these parameters as well as the corresponding volumes $w$. However, not all singularities are of interest in our problem since they might occur for parameters which are not located at the boundary of the region where local minima are possible. Moreover, once computed, it is not convenient to check if a singularity is of interest or not.\\
\\
Hence, for our numerical study, we proceed differently. Consider some parameters as fixed, for instance, like in the following, $t_3,\kappa,w_1,w_2$. Then, each couple of angles $(\alpha_1,\alpha_2)\in(-\pi,0)\times(0,\pi)$ (that is $(z_1,z_2)\in\Rbb^*_-\times\Rbb^*_+$) uniquely determines $t_1,t_2,z_3,y$ such that the system \eqref{linesys} is fulfilled and the computation of these parameters is straightforward, indeed
\begin{align}
y&=\sqrt[3]{\frac{z_2(z_2^2+3)-z_1(z_1^2+3)}{w_1+w_2}},\\
z_3&=2\sinh\left(\frac{1}{3}\arcsinh\left(\frac{w_2z_1(z_1^2+3)+w_1(z_2^2+3)}{2(w_1+w_2)}\right)\right),\\
t_1&=\frac{t_3\sin\varphi_1+\kappa ys_2}{\sin\varphi_3},\\
t_2&=\frac{t_3\sin\varphi_2-\kappa ys_1}{\sin\varphi_3}.
\end{align}
We have of course $y>0$ and $z_1<z_3<z_2$. Then, if $t_1>0$ and $t_2>0$ and if the configuration is a local minimum, the couple $(t_1,t_2)$ can be recorded as being in the domain of existence of local minima for the energy. Processing a lot of couples of angles allow to have a fairly precise knowledge of the set of couples $(t_1,t_2)$ such that the energy has a local minimum.\\
\\
As a first example, we consider a surface tension $t_3=1$ and for different line tensions $\kappa$, we observe which surface tensions $t_1,t_2$ give rise to an interior minimum.
\begin{figure}[H]\label{minima}
\centering
\includegraphics[scale=0.3]{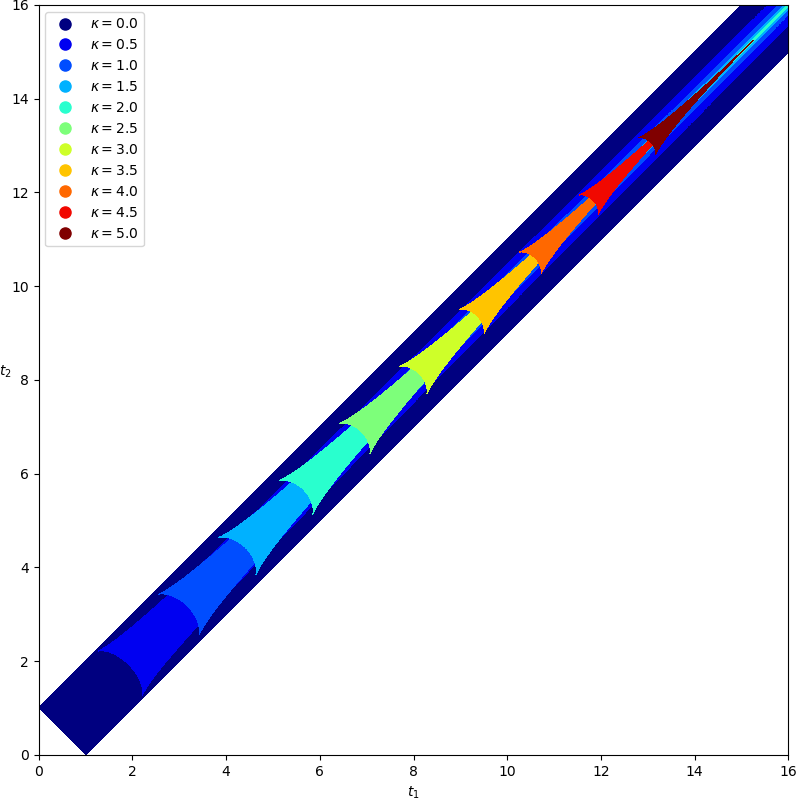}\includegraphics[scale=0.3]{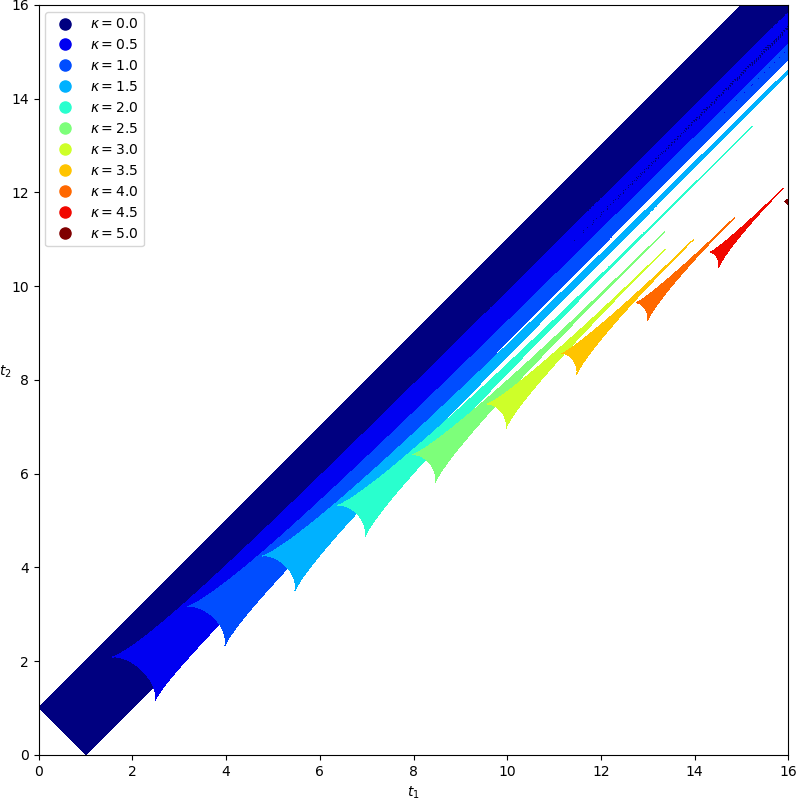}
\caption{Left: $w_1=0.5,w_2=0.5$. Right: $w_1=0.75,w_2=0.25$}
\end{figure}
From sections \ref{sysexindgen} and \ref{unique}, we already knew that for a doublet with only surface tensions, the domain of existence of a configuration does not depend on the volumes but only on the tensions (they have to fulfill the triangle inequalities). As expected, we can see that this is not true anymore for a doublet with line tension. Increasing the line tension also shrinks this domain. First unbounded, the domain becomes bounded when the line tension increases too much and then disappears when it is too high. When the cell volumes are not equal, we can also see that increasing the line tension slightly rotates the domain of existence of a configuration. We can observe these behaviors more into details and look at the value of the energy at the local minima. We still consider $t_3=1$ and we first consider equal volumes. For different values of the line tension $\kappa$, we again observe which surface tensions $t_1,t_2$ give rise to an interior minimum. For these local minima, we compare their energy to those of the three minima on the boundary. We denote by $E_k^0$, $k\in\{1,2,3\}$ the energy at $x_k=0$.\\
\\
Both on the domain of existence of a local minimum, the figure on the right shows the order of $E_1^0,E_2^0,E_3^0$ and the figure on the left compares the energy at the local minima to $E_1^0,E_2^0,E_3^0$. For comparison with a doublet without line tension, the black solid line represents the location of tensions where the triangle inequalities are violated. As a reminder, this is the same plot as in figure \ref{tensions_configurations}. Inside the domain $\{(t_1,t_2)\in\left(\Rbb^*_+\right)^2\,;\,|t_1-t_2|<t_3<t_1+t_2\}$, a cell doublet without line tension has a unique configuration of minimum energy with junction height different than $0$. Otherwise we either have internalization or separation.
\begin{figure}[H]
\centering
\includegraphics[scale=0.31]{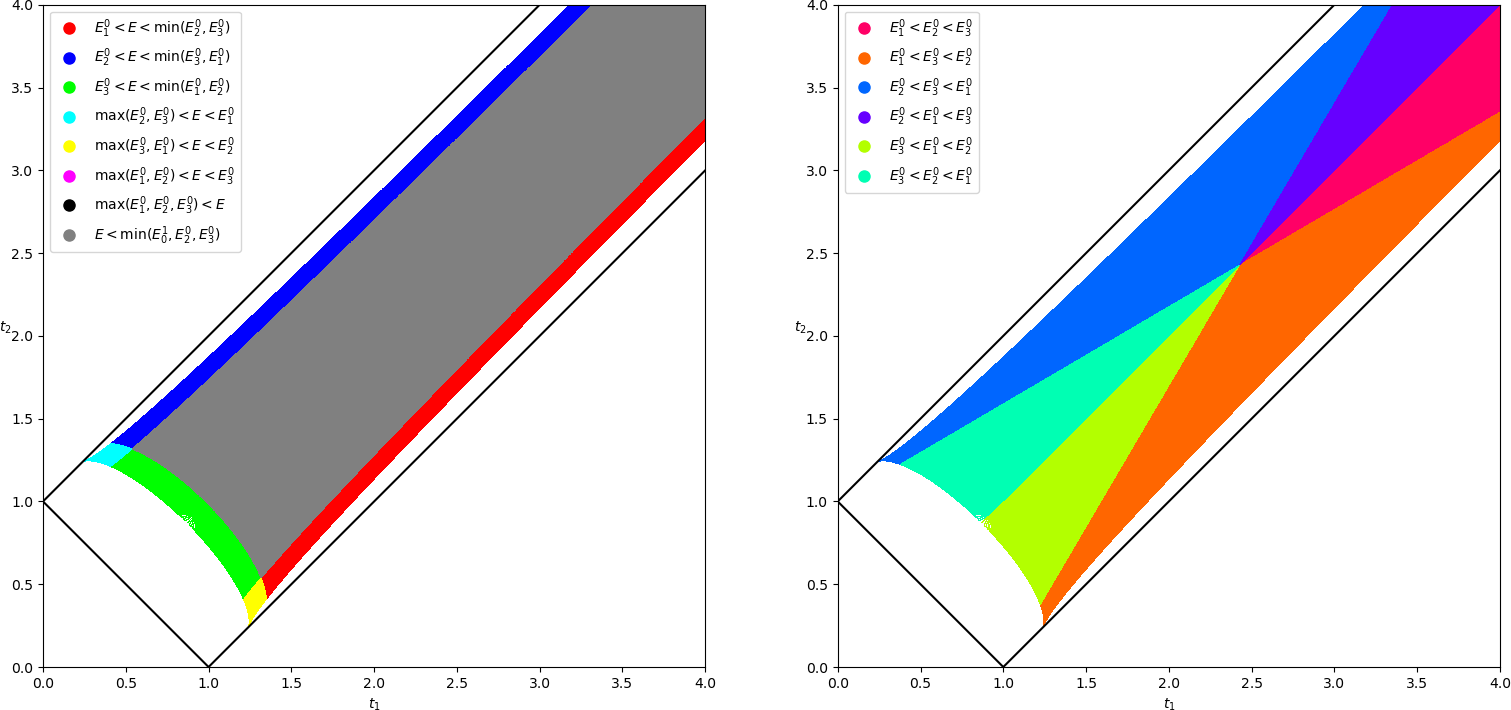}
\caption{$\kappa=0.1$, $t_3=1$, $w_1=w_2=0.5$.}
\end{figure}
Hence, in this case, we can see that most of local minima are also global ones. However, things might be different, increasing the line tension $\kappa$ makes the gray area become bounded and then decrease its size.
\begin{figure}[H]
\centering
\includegraphics[scale=0.31]{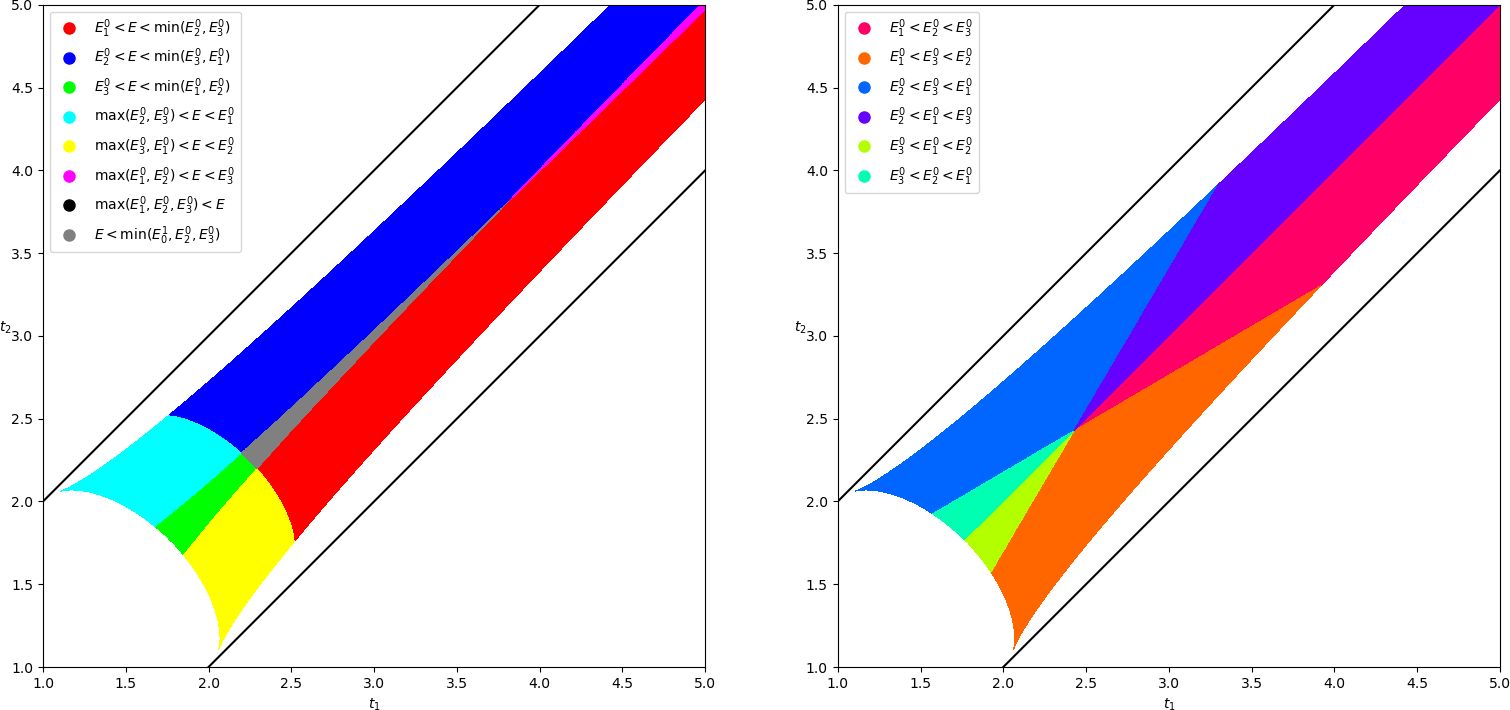}
\caption{$\kappa=0.44$, $t_3=1$, $w_1=w_2=0.5$.}
\end{figure}
It finally disappears so that in these cases, despite the existence of a local minimum such that $x_1,x_2,x_3\neq0$, the global minimum of the energy is always one of $E_1^0,E_2^0,E_3^0$.
\begin{figure}[H]
\centering
\includegraphics[scale=0.31]{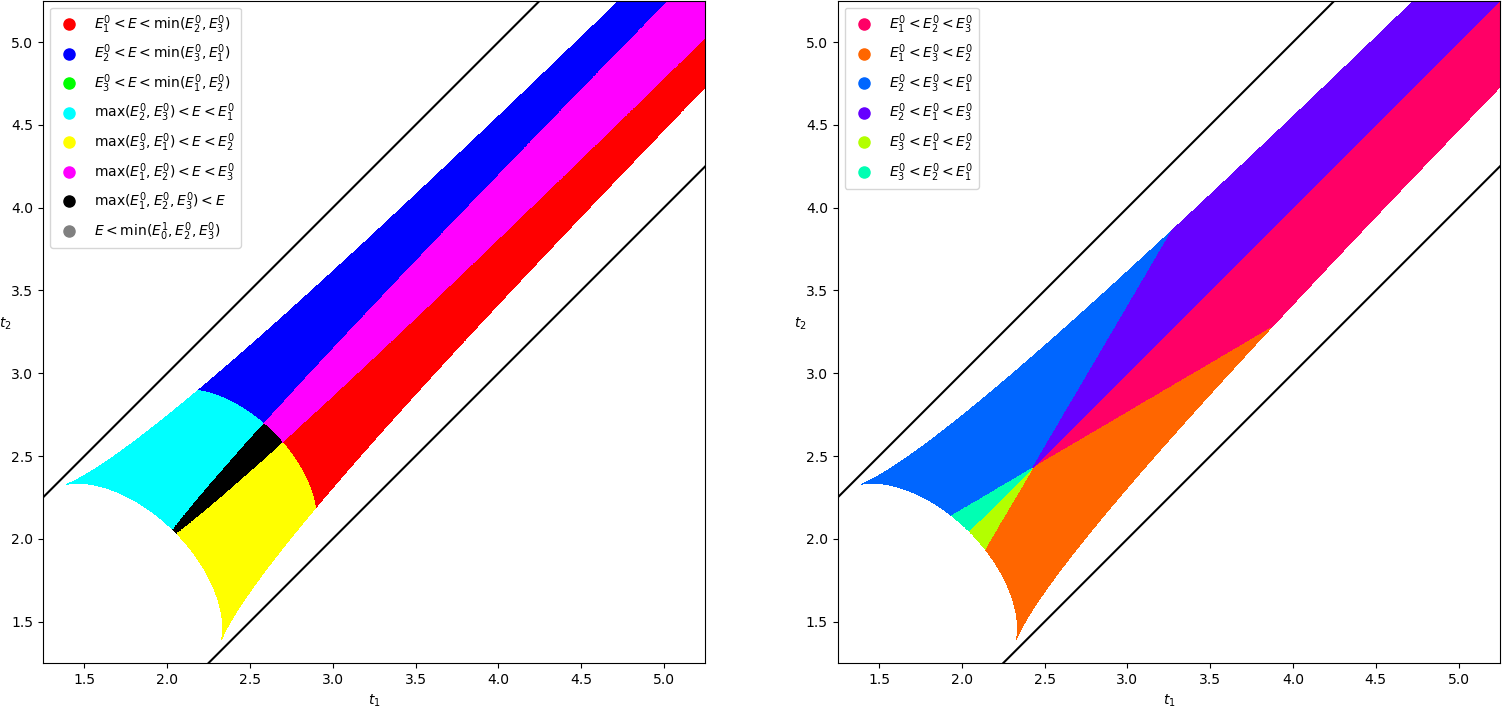}
\caption{$\kappa=0.55$, $t_3=1$, $w_1=w_2=0.5$.}
\end{figure}
On the previous figures, the domain of tensions $(t_1,t_2)$ where a local minimum exists is in fact unbounded so that $t_1,t_2$ can be chosen arbitrarily large if they do not differ too much. The same holds true for a cell doublet without line tension. However, in the presence of line tension, the domain becomes bounded if $\kappa$ increases too much. Precisely, in the case we are considering now ($w_1=w_2$), it happens when
\begin{equation}
\frac{\kappa}{\sqrt[3]{w_1+w_2}}>\frac{3t_3}{2}
\end{equation}
\begin{figure}[H]
\centering
\includegraphics[scale=0.31]{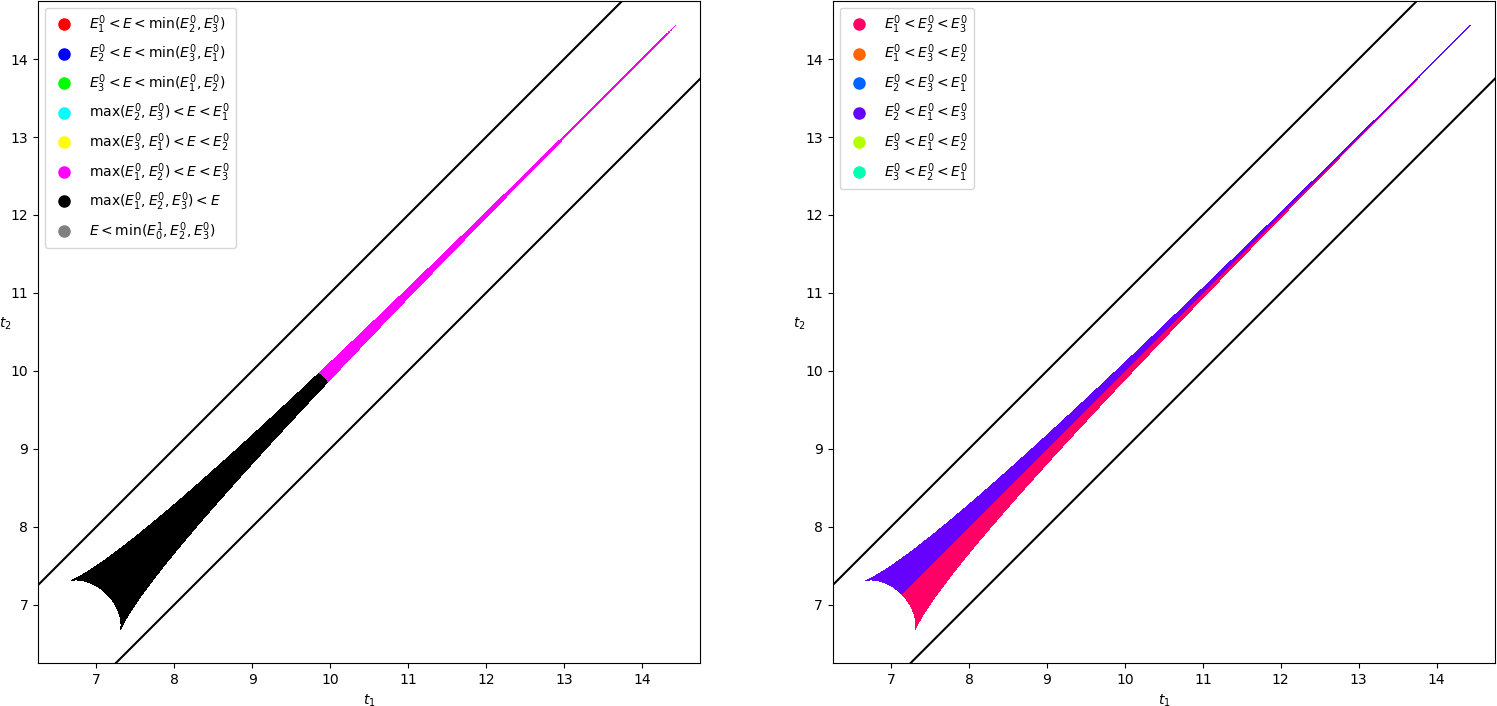}
\caption{$\kappa=2.6$, $t_3=1$, $w_1=w_2=0.5$.}
\end{figure}
If $\kappa$ increases more, then the energy of any local minimum is higher than the three energies for $h=0$.
\begin{figure}[H]
\centering
\includegraphics[scale=0.31]{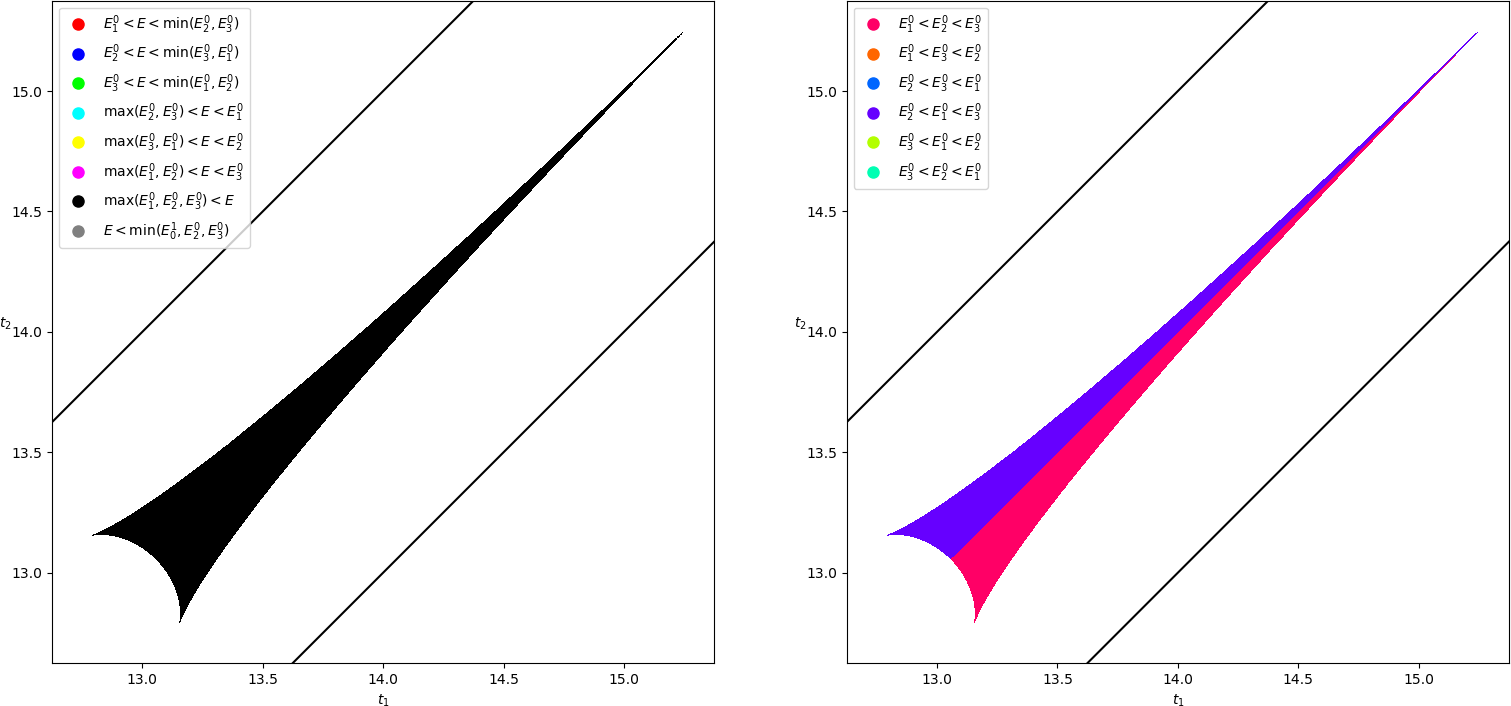}
\caption{$\kappa=5$, $t_3=1$, $w_1=w_2=0.5$.}
\end{figure}
As $\kappa$ increases, the set of couples $(t_1,t_2)$ such that the energy of a cell doublet has a local minimum is getting smaller and there is no more local minimum when $\kappa$ is too large, more precisely when
\begin{equation}\label{t3small}
\frac{\kappa}{\sqrt[3]{w_1+w_2}}>\frac{t_3}{4\sqrt[6]{14}}\left(43+5\sqrt{65}\right)\sqrt[6]{25\sqrt{65}-201},
\end{equation}
Where we still consider the case $w_1=w_2$. The common value $t$ of $t_1,t_2$ at the singularity (when equality holds in \eqref{t3small}) is given by
\begin{equation}
t=\frac{5t_3}{8}\left(25+3\sqrt{65}\right).
\end{equation}
In all previous cases, the existence of a local minimum is implying that the surface tensions fulfill the triangle inequalities. However, there is \textit{a priori} no reason for this fact to be true in general since we now have four forces ($\vec t_1,\vec t_2,\vec t_3,\vec\kappa$). Precisely, as already mentioned and seen on figure \ref{minima}, contrary to the case of a doublet without line tension where singularities are independent of the volumes, they play a role if a line tension exists. Hence, to illustrate what happens in this case, consider for instance $w_1=0.75,w_2=0.25$ and, as previously, a constant surface tension $t_3=1$ and the same values of the line tension $\kappa$. We then get the following plots (the scales are different and have been chosen to fit the region of local minima in the best way).
\begin{figure}[H]
\centering
\includegraphics[scale=0.31]{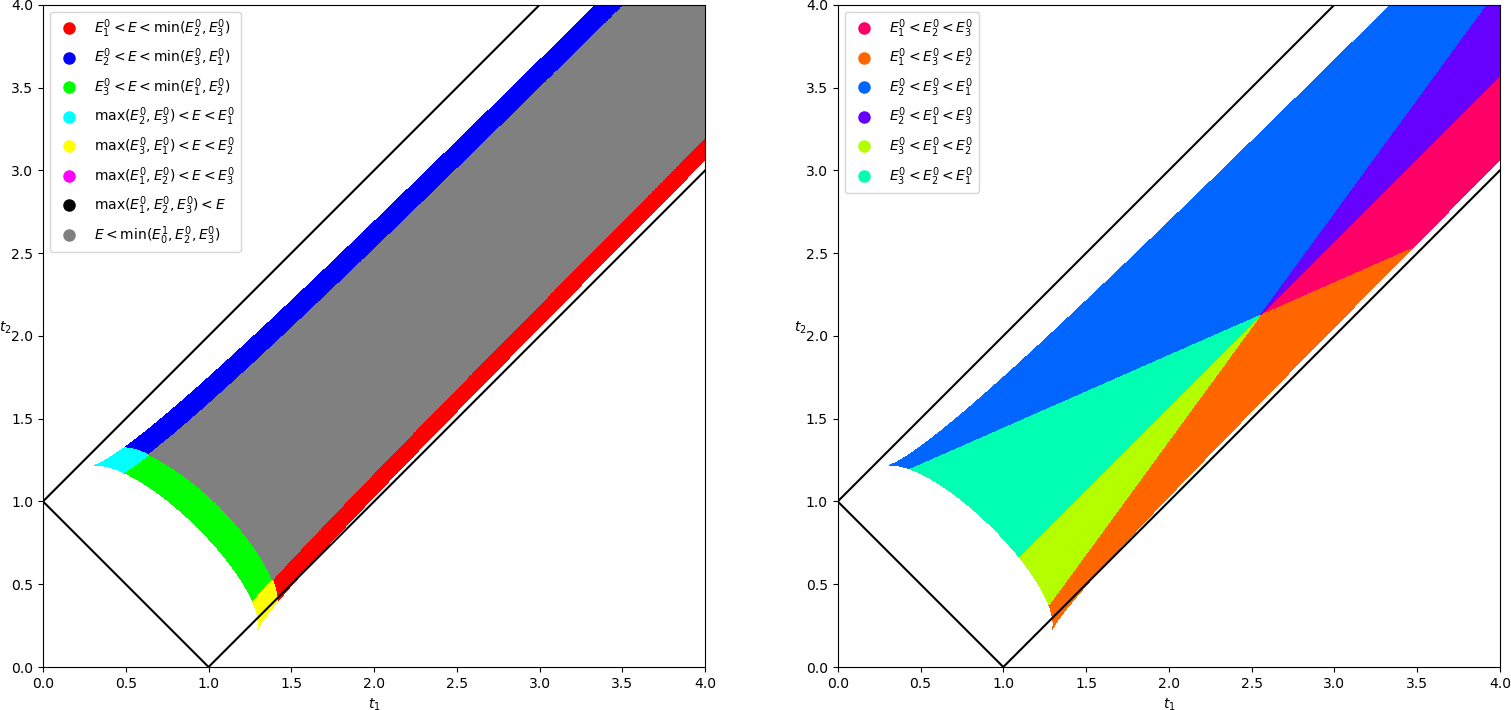}
\caption{$\kappa=0.1$, $t_3=1$, $w_1=0.75$, $w_2=0.25$.}
\end{figure}
We can already remark that for some local minima, the surface tensions do not fulfill the triangle inequalities (the region of local minima is slightly ``rotated''). This fact is even more obvious as the line tension increases.
\begin{figure}[H]
\centering
\includegraphics[scale=0.31]{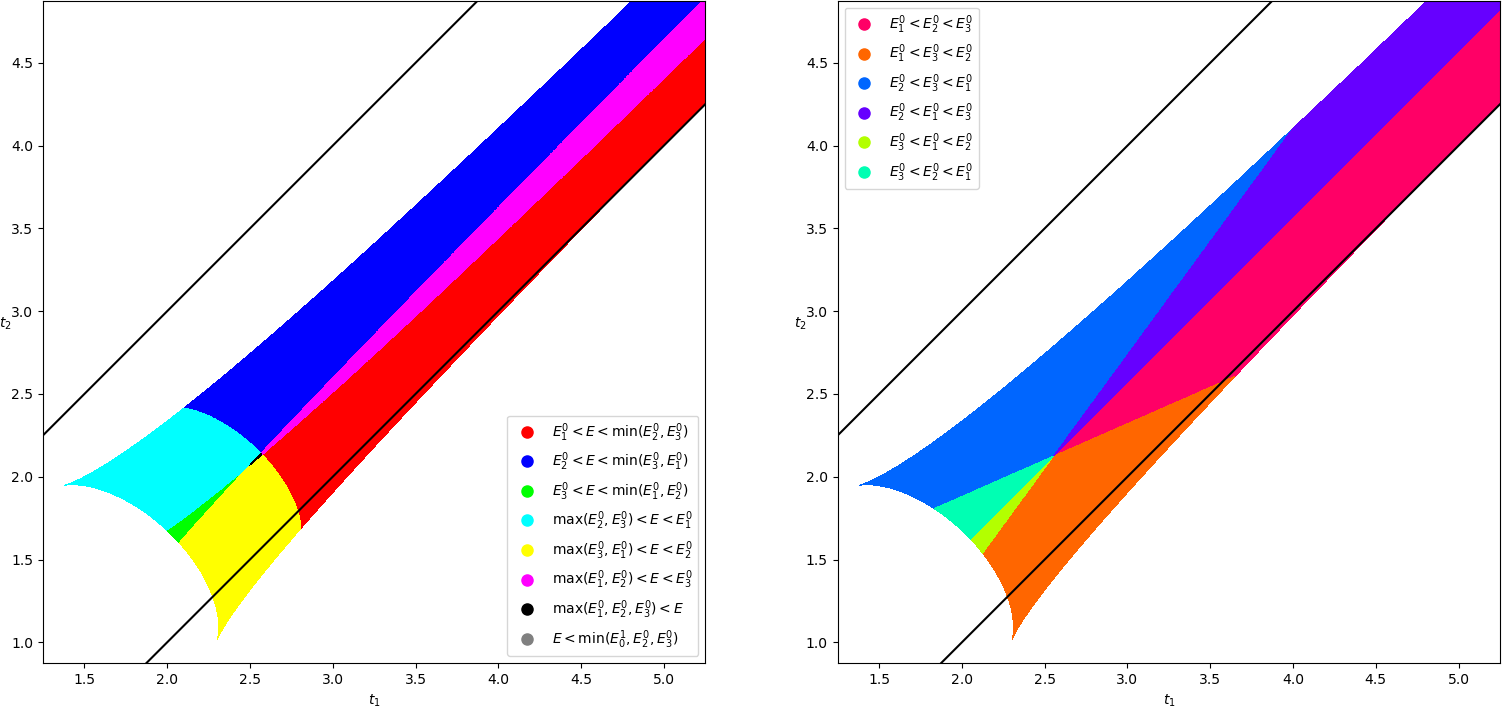}
\caption{$\kappa=0.44$, $t_3=1$, $w_1=0.75$, $w_2=0.25$.}
\end{figure}
\begin{figure}[H]
\centering
\includegraphics[scale=0.31]{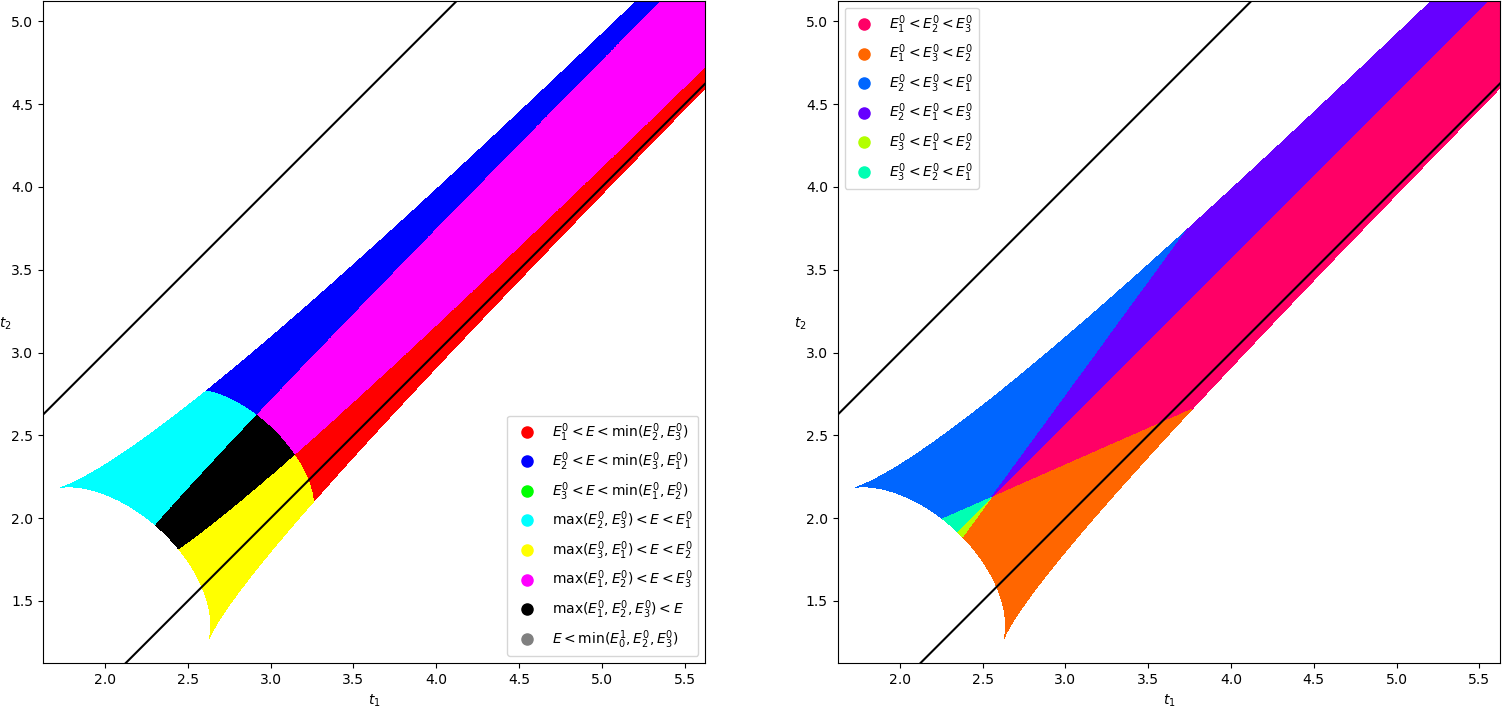}
\caption{$\kappa=0.55$, $t_3=1$, $w_1=0.75$, $w_2=0.25$.}
\end{figure}
For the two previous figures (``low'' line tension $\kappa$), besides the ``rotation'', we can also remark that the region of local minima is slightly narrower compared to the previous case of equal volumes. This fact becomes more visible when $\kappa$ increases. For $\kappa=2.6$, we get
\begin{figure}[H]
\centering
\includegraphics[scale=0.31]{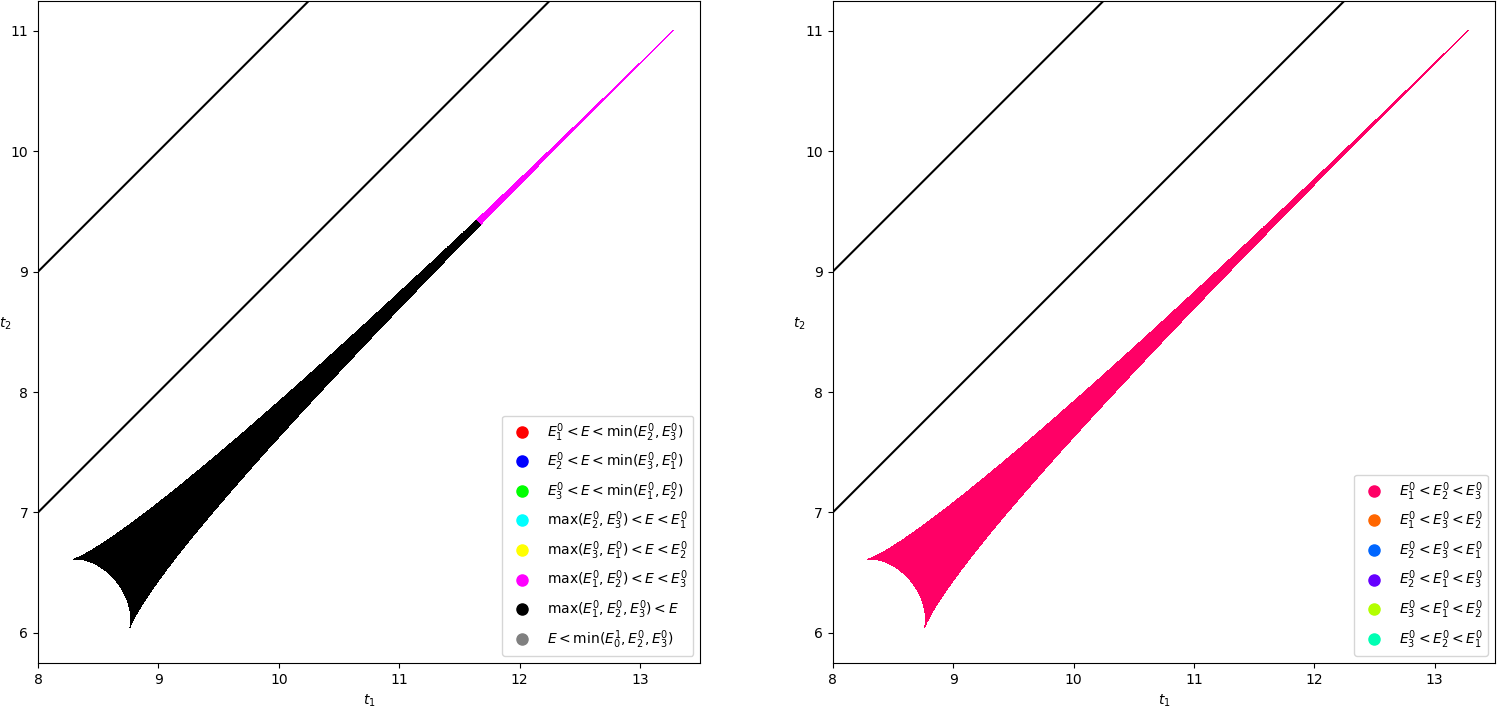}
\caption{$\kappa=2.6$, $t_3=1$, $w_1=0.75$, $w_2=0.25$.}
\end{figure}
Now the surface tensions in the region of local minima never fulfill the triangle inequalities. Moreover, the line tension for which the region becomes bounded has decreased compared to the equal volumes case. This value of $\kappa$ can be expressed as the root of a polynomial but it is probably not possible to express it in terms of radicals though. Now, considering $\kappa=5$, we can again remark that the region of local minima is much smaller.
\begin{figure}[H]
\centering
\includegraphics[scale=0.31]{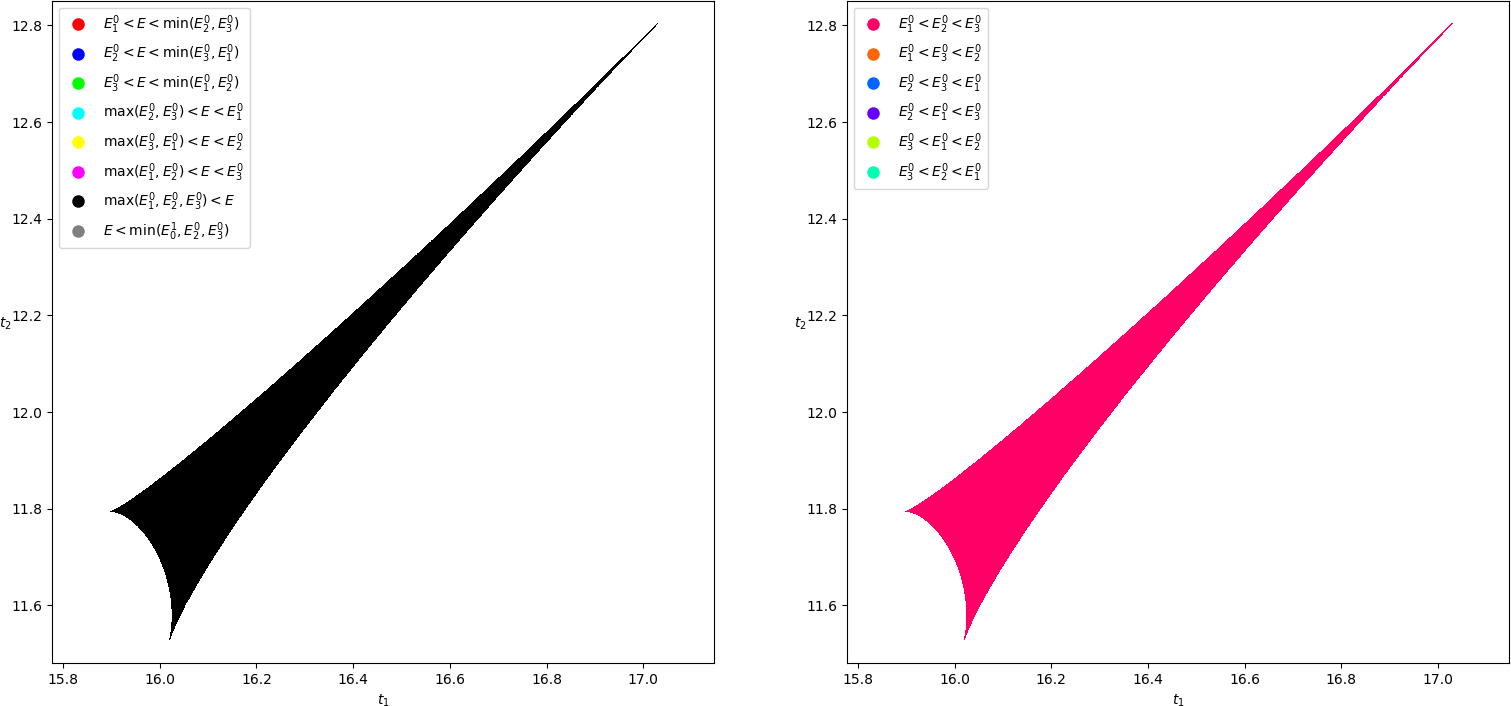}
\caption{$\kappa=5$, $t_3=1$, $w_1=0.75$, $w_2=0.25$.}
\end{figure}
As in the equal volumes case, when $\kappa$ increases enough, there is no more region of local minima. However, it happens for a lower line tension $\kappa$ than in the equal volumes case. As previously mentioned, it may not be possible to express the location of this singularity in terms of radicals but it can be written as a root of a polynomial.
\subsubsection{Bulging}\label{bulging}
For a cell doublet without line tension, we necessarily have $\sin\varphi_k>0$ for all $k\in\{1,2,3\}$. It is possible to have $\varphi_k\simeq\pi$ (but with $\varphi_k<\pi$) for a $k\in\{1,2,3\}$ if the surface tension $t_k$ is small compared to the others. For instance, with equal volumes, we get the following configurations when two of the tensions are twenty times larger than the last one.
\begin{figure}[H]
\centering
\includegraphics[scale=0.29]{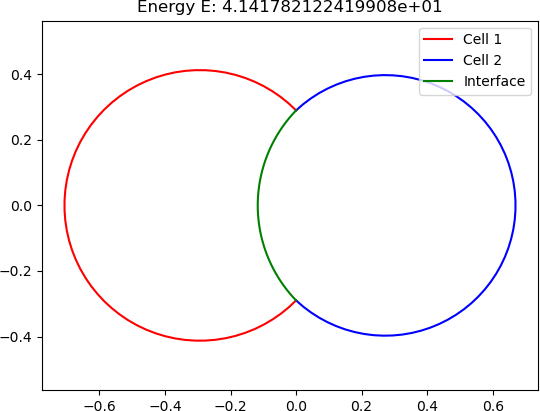}\includegraphics[scale=0.29]{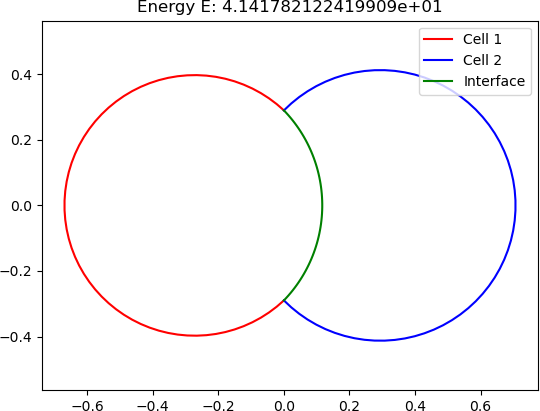}\includegraphics[scale=0.29]{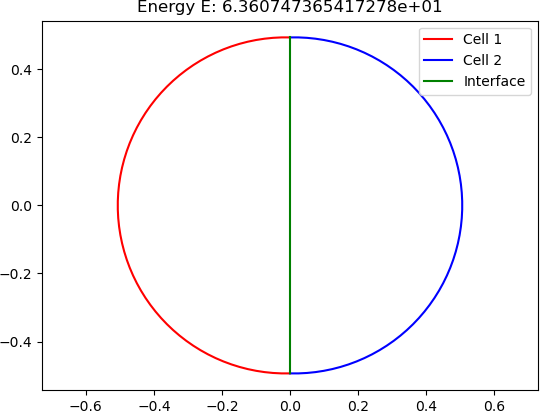}
\caption{From left to right:$\quad t_1\ll t_2,t_3\quad,\quad t_2\ll t_3,t_1\quad,\quad t_3\ll t_1,t_2$.}
\end{figure}
Now, this does not hold true anymore in the presence of line tension. Whereas we still have $\sin\varphi_3>0$, one of the other angles $\varphi_1,\varphi_2$ may be slightly larger then $\pi$. The cell doublet then looks like ``bulged''. On the following figures, we can see for which tensions bulging happens. Consider the same examples as previously, first with equal volumes.
\begin{figure}[H]
\centering
\includegraphics[scale=0.36]{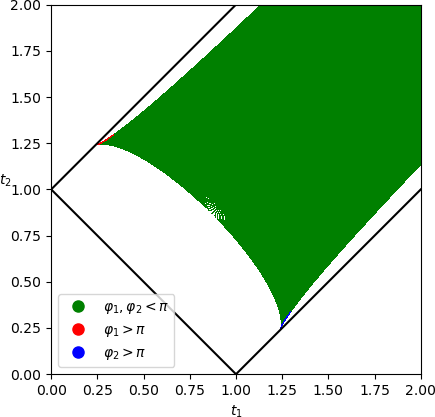}\includegraphics[scale=0.36]{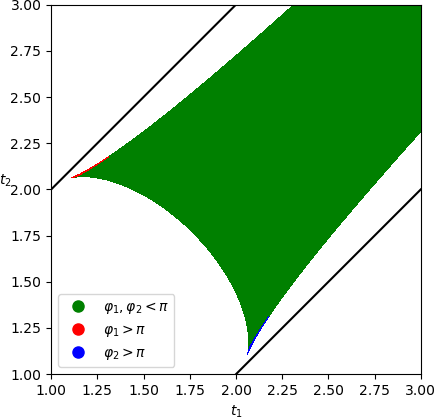}\includegraphics[scale=0.36]{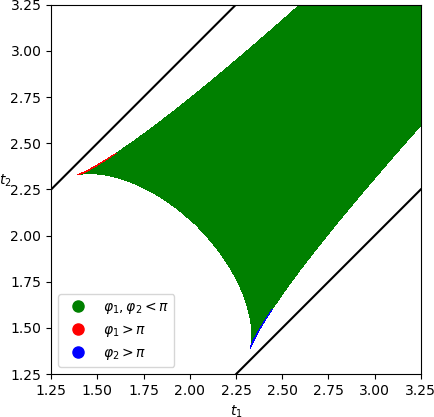}
\caption{$t_3=1$ and from left to right:$\quad\kappa=0.1\quad,\quad\kappa=0.44\quad,\quad\kappa=0.55$.}
\end{figure}
The area of bulging is rather small and located close to the boundary of the domain. For larger $\kappa$ there is no more possible bulging. Since the region of bulging is very small and located close to the singularities, one can wonder if the previous results are not due to numerical errors. In order to be sure that this is not the case, one can check that the region of local minima can contain solutions such that $\sin\varphi_1=0$ or $\sin\varphi_2=0$ so that a slight modification of the tensions can lead to $\sin\varphi_1<0$ or $\sin\varphi_2<0$. To do so, the system \eqref{linesys} can easily be solved exactly when assuming for instance $\sin\varphi_1=0$. Indeed, we then have $z_3=-1/z_2$ and the system writes
\begin{align}
t_1\frac{1-z_1^2}{1+z_1^2}+(t_2-t_3)\frac{1-z_2^2}{1+z_2^2}+\kappa y&=0,\\
t_1\frac{z_1}{1+z_1^2}+(t_2-t_3)\frac{z_2}{1+z_2^2}&=0,\\
-\frac{1+3z_2^2}{z_2^3}-z_1(z_1^2+3)&=w_1y^3,\\
\left(\frac{1+z_2^2}{z_2}\right)^3&=w_2y^3.
\end{align}
Hence we immediately get
\begin{equation}
y=\frac{1+z_2^2}{z_2\sqrt[3]{w_2}}.
\end{equation}
Replacing in the other equations, the first two can be used to find $z_1,z_2$ and the last one gives a condition on $t_1,t_2,t_3,w_1,w_2$. We have
\begin{equation}
z_2=\frac{1}{2\kappa}\left((t_2-t_3)\sqrt[3]{w_2}+r\pm\sqrt{\sqrt[3]{w_2}((t_2-t_3)^2+t_1^2)+2r(t_2-t_3)}\right),
\end{equation}
With
\begin{equation}
r=\pm\sqrt{(t_1\sqrt[3]{w_2}+2\kappa)(t_1\sqrt[3]{w_2}-2\kappa)}.
\end{equation}
And
\begin{equation}
z_1=-\frac{2\kappa(t_2-t_3)z_2}{\kappa(t_1-t_2+t_3)z_2^2-\sqrt[3]{w_2}(t_1^2-(t_2-t_3)^2)z_2+\kappa(t_1+t_2-t_3)}.
\end{equation}
For equal volumes, $\kappa=0.1$, $t_3=1$ and $t_2=5/4$ (where bulging seems to happen), $t_1$ is then the largest real positive root of $P$ where
\begin{equation}
\begin{aligned}
P&=X^{10}-\frac{1}{8}X^8-\frac{3}{64}X^7+\left(\frac{1}{256}+\frac{3\sqrt[3]{2}}{10000}+\frac{\sqrt[3]{4}}{100}\right)X^6+\left(\frac{3}{512}+\frac{27\sqrt[3]{2}}{40000}\right)X^5+\\
&\left(\frac{15881}{32000000}+\frac{3\sqrt[3]{2}}{40000}\right)X^4+\left(-\frac{3}{16384}+\frac{27\sqrt[3]{2}}{640000}-\frac{3\sqrt[3]{4}}{6400}\right)X^3+\\
&\left(-\frac{15497}{256000000}-\frac{33\sqrt[3]{2}}{2560000}\right)X^2+\frac{15881}{8192000000}-\frac{9\sqrt[3]{2}}{10240000}+\frac{15629\sqrt[3]{4}}{3200000000}.
\end{aligned}
\end{equation}
Using a Computer Algebra System, the different quantities can then be evaluated with arbitrary high precision. We get $t_1=0.271244499897851$, $z_1=-2.031331771464625$, $z_2=1.756721787151541$, $z_3=-0.569242100436099$, $y=2.930530863266979$. The trace of the Hessian matrix is given by\\ $26.158972256436426$ and it determinant by $4.615158111925040$. The configuration is then a local minimizer. By slightly decreasing $t_1$ ($-0.003$ maximum) or slightly increasing $t_2$ ($+0.002$ maximum), one then obtains a configuration such that $\varphi_1>\pi$ which is a local minimizer. If we now look at what happens for different volumes ($w_1=0.75$, $w_2=0.25$), we have similar results.
\begin{figure}[H]
\centering
\includegraphics[scale=0.36]{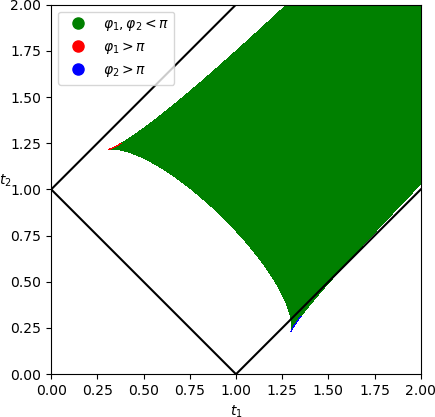}\includegraphics[scale=0.36]{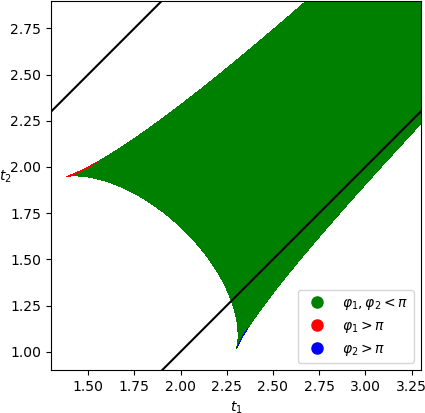}\includegraphics[scale=0.36]{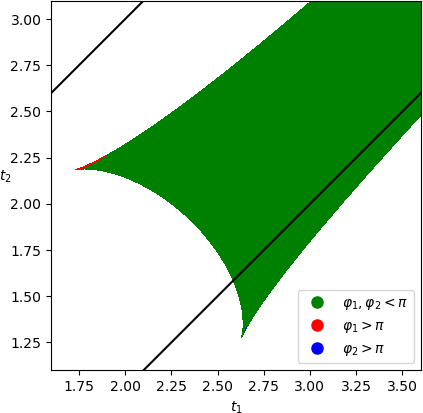}
\caption{$t_3=1$ and from left to right:$\quad\kappa=0.1\quad,\quad\kappa=0.44\quad,\quad\kappa=0.55$.}
\end{figure}
For larger $\kappa$, the area $\varphi_2>\pi$ disappears first (it has almost disappeared for $\kappa=0.55$), then the area $\varphi_1>\pi$ disappears as well.\\
\\
Running several computations with a large number of parameters, it seems that bulging is limited to about $2.6^\circ$ above flat angle. More important bulging which can be observed during some experiments may fit with the model though. Indeed, the case studied here is only a double cell and some things way differ when observing an early embryo for instance. We can also remark that the constraints (volume constraints here) play an important role in the configurations which can be observed. Different constraints may possibly lead to more important bulging. To end up with these remarks, bulging remarked in experiments can be due to other physical phenomenons and/or combined effects of line tension studied here with something else. On the next figure, we can see how a double cell with maximum bulging appears. Precisely, for the following parameters
\begin{equation*}
t_1=3.544814028\quad,\quad t_2=3.838944848\quad,\quad t_3=1.730430441,
\end{equation*}
\begin{equation*}
\kappa=1.000000000\quad,\quad w_1=0.820008308\quad,\quad w_2=0.179991692,
\end{equation*}
We have $\varphi_1=182.590653^\circ$. The double cell has the following configuration
\begin{figure}[H]
\centering
\includegraphics[scale=0.45]{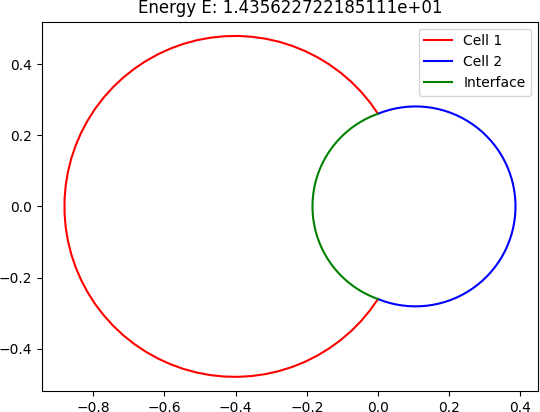}
\caption{Double cell with bulging.}\label{bulgingconf}
\end{figure}
\subsubsection{Tension inference}
As we have seen in section \ref{teninf}, for a doublet with only surface tensions, it is possible to infer the tensions (up to a common scaling factor). Moreover, from the previous section \ref{bulging}, line tensions have effects which make a configuration possibly differ from one without line tension. Thus, one can wonder if inferring the tensions is possible in this more general case. In fact it seems to be an impossible task in general. Indeed, consider a doublet with parameters $t_1,t_2,t_3,\kappa,w_1,w_2$ at local minimum. Of course, as in the case of double cells with only surface tensions, common scaling of the tensions does not change the configuration, but we have more since we have now four tensions and still only two equations involving the tensions which have to be fulfilled by the configuration parameters, namely, for a given configuration at local minimum, we have
\begin{align*}
t_1c_1+t_2c_2+t_3c_3+\kappa y&=0,\\
t_1s_1+t_2s_2+t_3s_3&=0.
\end{align*}
This means that we have a vector plane (space of dimension $2$ in the space of dimension $4$) of tensions $(\hat t_1,\hat t_2,\hat t_3,\hat \kappa)$ fulfilling these two equations. More precisely, we have
\begin{align*}
\hat t_1&=\lambda t_1+\mu\sin\varphi_1,\\
\hat t_2&=\lambda t_2+\mu\sin\varphi_2,\\
\hat t_3&=\lambda t_3+\mu\sin\varphi_3,\\
\hat \kappa&=\lambda \kappa,
\end{align*}
With $\mu\in\Rbb$ and $\dsp\lambda\in\Rbb_+^*$ such that $\dsp \lambda>-\mu\min_{k\in\{1,2,3\}}\sin\varphi_k/t_k$. $\lambda,\mu$ have also to be such that the configuration is not only a critical point but also a local minimum, a necessary condition being that \eqref{trhessht} and \eqref{dethessht} are non positive. If the inequalities are strict for $t_1,t_2,t_3$, then there always exists a neighborhood of $\mathcal{V}$ of $0$ such that $\hat t_1>0,\hat t_2>0,\hat t_3>0,\hat \kappa>0$ and the corresponding configuration is a local minimum for $\lambda\in\Rbb_+^*$ and $\mu\in\lambda\mathcal{V}$.\\
\\
The specially interesting case of this non uniqueness is when we take $\hat \kappa=0$ (that is $\lambda=0$), then the solutions $(\hat t_1,\hat t_2,\hat t_3)$ are simply the tensions verifying Lami's theorem \eqref{lami}
\begin{equation}
\frac{\hat t_1}{\sin\varphi_1}=\frac{\hat t_2}{\sin\varphi_2}=\frac{\hat t_3}{\sin\varphi_3}.
\end{equation}
Consider for instance the following example

\begin{minipage}{6cm}
\begin{figure}[H]
\centering
\includegraphics[width=6cm]{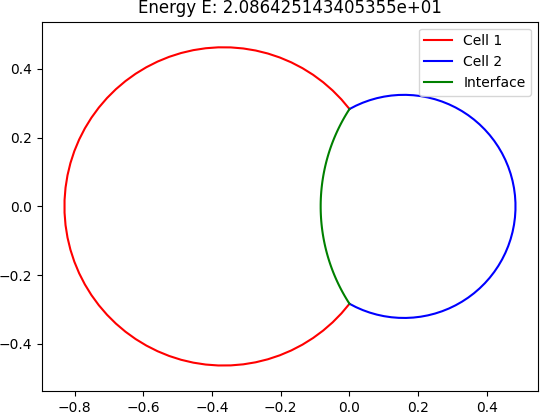}
\caption{$t_1=5$, $t_2=6$, $t_3=4$, $\kappa=1$, $w_1=0.75$, $w_2=0.25$.}\label{cellwithk}
\end{figure}
\end{minipage}
\begin{minipage}{6cm}
\centerline{Angles are given by}
\begin{align}
\varphi_1&=52.2^\circ,\\
\varphi_2&=109.3^\circ,\\
\varphi_3&=98.5^\circ.
\end{align}
\end{minipage}
\\
\\
\\
Then the double cell with surface tensions $t_1=4.2$, $t_2=8.5$, $t_3=8.9$ and with the same volumes $w_1=0.75$, $w_2=0.25$ has exactly the same configuration as the previous doublet with line tension.
\begin{figure}[H]
\centering
\includegraphics[width=6cm]{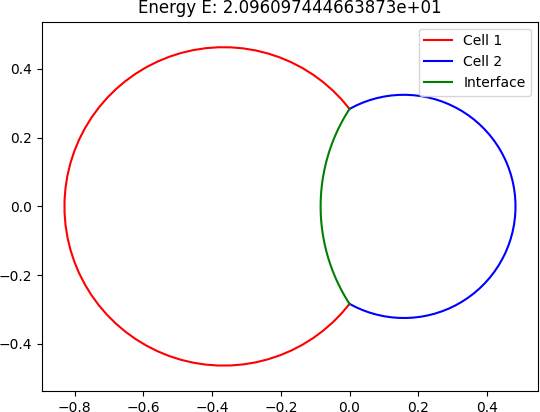}
\caption{$t_1=4.2$, $t_2=8.5$, $t_3=8.9$, $w_1=0.75$, $w_2=0.25$.}\label{cellwithoutk}
\end{figure}
So, though line tensions can produce typical bulging configurations, they mainly produce classical configurations which are geometrically indistinguishable from a doublet without line tension (the pressures differ though).
\bibliography{doublet}

\begin{thebibliography}{10}

\bibitem{abramowitz_stegun}
Milton Abramowitz and Irene~A. Stegun.
\newblock {\em Handbook of mathematical functions with formulas, graphs, and
  mathematical tables}, volume~55 of {\em National Bureau of Standards Applied
  Mathematics Series}.
\newblock For sale by the Superintendent of Documents, U.S. Government Printing
  Office, Washington, D.C., 1964.

\bibitem{barrat}
Jean-Louis Barrat and J.-P. Hansen.
\newblock Basic concepts for simple and complex liquids.
\newblock {\em Physics Today}, 58:56--57, 01 2005.

\bibitem{brioschi}
F.~Brioschi.
\newblock Sulla risoluzione delle equazioni del quinto grado.
\newblock {\em Annali di Matematica Pura ed Applicata}, 1858.

\bibitem{chaikin_lubensky}
P.~M. Chaikin and T.~C. Lubensky.
\newblock {\em Principles of Condensed Matter Physics}.
\newblock Cambridge University Press, Cambridge, 1995.

\bibitem{cox_sturmfels_manocha_sederberg_kramer_laubenbaches_thomas}
David Cox, Bernd Sturmfels, Dinesh Manocha, Thomas Sederberg, Z.~Kramer,
  R.~Laubenbaches, and R.~Thomas.
\newblock {\em Applications of Computational Algebraic Geometry}.
\newblock 03 2015.

\bibitem{cox}
David~A. Cox.
\newblock {\em Galois theory [2nd ed.]}.
\newblock Wiley, New Jersey, 2nd ed. edition, 2012.

\bibitem{gennes_brochard-wyart_quere}
P.-G. de~Gennes, F.~Brochard-Wyart, and D.~Qu{\'e}r{\'e}.
\newblock {\em Gouttes, bulles, perles et ondes}.
\newblock {\'E}chelles. Belin, 2002.

\bibitem{singular}
Wolfram Decker, Gert-Martin Greuel, Gerhard Pfister, and Hans Sch\"onemann.
\newblock {\sc Singular} {4-3-0} --- {A} computer algebra system for polynomial
  computations.
\newblock \url{http://www.singular.uni-kl.de}, 2022.

\bibitem{hermite}
C.~Hermite.
\newblock Sur la r{\'e}solution de l'{\'e}quation du cinqui{\`e}me degr{\'e}.
  {I}n {{\OE}}uvres de {C}harles {H}ermite.
\newblock pages 5--12, 2009.

\bibitem{hutchings_morgan_ritore_ros}
Michael Hutchings, Frank Morgan, Manuel Ritoré, and Antonio Ros.
\newblock Proof of the double bubble conjecture.
\newblock {\em Annals of Mathematics}, 155(2):459--489, 2002.

\bibitem{israelachvili}
J.~N. Israelachvili.
\newblock {\em Intermolecular and Surface Forces, 3rd Edition}.
\newblock Elsevier Academic Press Inc, 525 B Street, Suite 1900, San Diego, CA
  92101-4495 USA, 2011.

\bibitem{klein_morrice}
Felix~Christian Klein and George~Gavin Morrice.
\newblock Lectures on the icosahedron and the solution of equations of the
  fifth degree.
\newblock 2003.

\bibitem{lawlor}
Gary~R Lawlor.
\newblock Double bubbles for immiscible fluids in $\mathbb{R}^n$.
\newblock {\em Journal of Geometric Analysis}, 24:190--204, 2014.

\bibitem{michalek_sturmfels}
Mateusz Michałek and Bernd Sturmfels.
\newblock {\em Invitation to nonlinear algebra}.
\newblock Number 211 in Graduate studies in mathematics. American Mathematical
  Society, Providence, Rhode Island, 2021.

\bibitem{nash}
Oliver Nash.
\newblock On {K}lein’s icosahedral solution of the quintic.
\newblock {\em Expositiones Mathematicae}, 32(2):99--120, 2014.

\bibitem{safran}
Samuel~A. Safran.
\newblock Statistical thermodynamics of surfaces, interfaces, and membranes.
\newblock 1994.

\bibitem{sturmfels_1995}
Bernd Sturmfels.
\newblock Gr{\"o}bner bases and convex polytopes.
\newblock 1995.

\bibitem{sturmfels_2002}
Bernd Sturmfels.
\newblock Solving systems of polynomial equations.
\newblock In {\em American Mathematical Society, CBMS Regional Conferences
  series, no 97}, 2002.

\bibitem{tschirnhaus_green}
Ehrenfried Tschirnhaus and R.~Green.
\newblock A method for removing all intermediate terms from a given equation.
\newblock {\em ACM Sigsam Bulletin}, 37:1--3, 03 2003.

\bibitem{verschelde}
Jan Verschelde.
\newblock Algorithm 795: {PHC}pack: A general-purpose solver for polynomial
  systems by homotopy continuation.
\newblock {\em ACM Trans. Math. Softw.}, 25(2):251–276, jun 1999.

\end{thebibliography}
\end{document}